HIGHER SCHOOL OF ECONOMICS
NATIONAL RESEARCH UNIVERSITY

*Sergey Shvydun*

# NORMATIVE PROPERTIES OF MULTI-CRITERIA CHOICE PROCEDURES AND THEIR SUPERPOSITIONS: II

Two-stage superposition choice procedures, which sequentially apply two choice procedures so that the result of the first choice procedure is the input for the second choice procedure, are studied. We define which of them satisfy given normative conditions, showing how a final choice is changed due to the changes of preferences or a set of feasible alternatives. A theorem is proved showing which normative conditions are satisfied for two-stage choice procedures based on different scoring rules, rules, using majority relation, value function and tournament matrix. A complexity of two-stage choice procedures as well as its runtime on real data are evaluated.

*Shvydun Sergey* – National Research University Higher School of Economics, Moscow, Russia; V.A. Trapeznikov Institute of Control Sciences of Russian Academy of Sciences, Moscow, Russia.

The work was partially financed by the International Laboratory of Decision Choice and Analysis (DeCAn Lab) of the National Research University Higher School of Economics and by V.A. Trapeznikov Institute of Control Sciences of Russian Academy of Sciences. The research was partially supported also by the Russian Foundation for Basic Research (RFBR) (grant № 12-00-00226 “Choice models based on superposition”).

The author would like to thank Professor F. Aleskerov for the formulation of the research problem. The author also thank Professor V. Volsky for his comments and suggestions, which helped the author to improve the paper.

## Introduction

The choice of the best alternatives among a set of all possible alternatives has been a matter for study, analysis and debates for a long time. It is hardly to find any sphere where this problem did not occur.

There are a lot of different choice procedures that allow to choose and rank alternatives from the initial set [1-10]. In this part, we consider choice procedures of a special type based on the superposition principle. Let us remind that by superposition of two choice functions $C_1(\cdot)$ and $C_2(\cdot)$ we mean a binary operation $\odot$, the result of which is a new function $C^*(\cdot)=C_2(\cdot)\odot C_1(\cdot)$, having the form $\forall X \epsilon 2^A$ $C^*(X)=C_2(C_1(X))$ [1]. In other words, superposition consists in sequential application of choice functions where the result of the previous choice function $C_1$ is the input for the next choice function $C_2$. It is necessary to mention that the change of the order of functions may lead to completely different results, as the superposition operation is not commutative.

The interest in superposition of choice procedures can be explained by several reasons. First, most existing accurate choice procedures have a high computational complexity so they cannot be applied in the cases when the number of alternatives or/and criteria is very large. The use of superposition allows to reduce the complexity by applying choice procedures with a low computational complexity on first stages and more accurate choice procedures on final stages. Thus, the results can be obtained in a reasonable time. Second, there are a lot of situations when after applying some choice procedures the remaining set of alternatives is too large. The use of superposition allows to avoid such situations through the use of additional choice procedures.

The change of presentation, a set of criteria or criterial values of some alternatives can affect the final choice. Consequently, there is a need for more detailed study of existing choice procedures and for understanding which of them can be used in a particular case.

Thus, we consider the two-stage superposition choice procedures based on scoring rules, rules, using majority relation, value function and tournament matrix. The main focus of the paper is the study of its properties, its computational complexity and its runtime on real data. The study of the properties of two-stage superposition choice procedures is based on the study of the properties of different multi-criteria choice procedures which is done in the first part of the study [11].

This part is divided into several sections. First, some background information on two-stage choice procedures is given. Second, we form a list of two-stage choice procedures Then, we study the properties of two-stage superposition choice procedures. Finally, a computational and run-time complexity of studied procedures is given.

## A survey of the literature

In [1] two-stage extremizational choice procedures that consist of scalar or vector choice procedures were studied. The description of such choice procedures is given below.

*Definition 1.* A choice function $C(\cdot)$ is called rationalized by scalar criterion $\varphi$ (or simply scalar), if $\forall X \subseteq A$

$$C(X) = \{y \in X | \nexists x \in X: \varphi(x) > \varphi(y)\}.$$

Suppose now that all alternatives from $A$ are mapped into $n$ criterial scales instead of one criterial scale. Denote by $\vec{\varphi}$ a set of $n$ criterial scales $\vec{\varphi} = (\varphi_1, \ldots, \varphi_n)$, where $\vec{\varphi}$ is a «vector criterion», and $\varphi_i$, where $i = 1, 2, \ldots, n$, - its component.

*Definition 2.* A choice function $C(\cdot)$ is called rationalized by vector criterion $\vec{\varphi} = (\varphi_1, \ldots, \varphi_n)$ or simply vector, if $\forall X \subseteq A$

$$C(X) = \{y \in X | \nexists x \in X: \vec{\varphi}(x) > \vec{\varphi}(y)\},$$

where

$$\vec{\varphi}(x) > \vec{\varphi}(y) \Leftrightarrow \forall i \in \{1, \ldots, n\}\ \varphi_i(y) > \varphi_i(x)$$

Thus, there are 4 main types of two-stage extremizational choice procedures

1. scalar-scalar choice procedure;
2. scalar-vector choice procedure;
3. vector-scalar choice procedure;
4. vector-vector choice procedure.

As the single–criterion extremizational choice procedure is a special case of the vector choice procedure, it is clear that the choice procedure 1 is a special case of choice procedures 2 and 3, which are special cases of the procedure 4. In [1] it was also shown that scalar-scalar choice procedure (type 1) is equivalently reducible to the usual single-criterion extremizational choice procedure and two-stage scalar-vector choice procedure - to the usual multi-criteria extremizational choice procedure. To define under what circumstances the procedure of type 2 is equivalent to the scalar procedure, a notion of $\psi$-triad was introduced in [1].

*Definition 3*. Let *u, v, w*$\in$*A*, where $A$ is a set of alternatives. These alternatives are said to form a $\psi$-triad if $\psi(u) < \psi(v)$, $\psi(u)\ \chi\ \psi(w)$ and $\psi(v)\ \chi\ \psi(w)$, where the inequality is understood as some vector (component-wise) inequality, and $\chi$ stands for the independent relation introduced as

$$\psi(u) \chi\ \psi(w) \Leftrightarrow \exists i, j \in I = \{1, \ldots, n\}: \psi_i(u) \geq \psi_i(w)\ and\ \psi_j(u) \leq \psi_j(w),$$

where $I -$ a set of criteria, $\psi$ – vector criterion.

The scalar-vector two-stage choice procedure (type 2) is equivalent to the one-stage single-criterion extremizational choice procedure if and only if the set $X$ lacks $\psi$-triads for any $X \in 2^A$ and the criterion $\varphi$ has the same value over all alternatives from this set, that is $\forall x \in X$ $\varphi(x) = const$ [1].

As for the vector-vector choice procedure of type 4, it was shown in [1] that this procedure is far from being always reducible to the usual multi-criteria (all the more so, to single-criterion) choice procedure. Reducibility is possible only for some special relative positions of the alternatives in criterial spaces $\varphi$ and $\psi$.

In [1] there were also defined under which conditions choice procedures of type 3 and 4 can be reduced to pair-dominant choice procedures. Let us remind that a choice procedure $C(\cdot)$ is called pair-dominant or rationalized by binary relation $P$ if $\forall X \subseteq A$

$$C(X) = \{y \in X | \nexists x \in X: xPy\} \qquad (1)$$

In other words, pair-dominant choice procedure is a procedure for which the rational choice consists in the choice of undominated by their pairwise comparisons alternatives.

To define under which conditions choice procedures of type 4 can be reduced to pair-dominant choice procedure, a notion of uncoordinated triad was introduced.

*Definition 4.* Let *u,v,w*$\in$*A* make up a $\varphi$-triad referred to as a $\varphi$, $\psi$-uncoordinated triad of the first, or second, or third type if correspondingly

1. $\psi(w) < \psi(u)$, $\psi(w) \not< \psi(v)$ or
2. $\psi(u)\,\chi\,\psi(w)$, $\psi(v) < \psi(w)$ or
3. $\psi(u)\,\chi\,\psi(w)$, $\psi(v)\,\chi\,\psi(w)$,

where $\varphi$, $\psi$ – vector criteria, $\chi$ – independed relation, $\not<$ is a violation of the vector inequality <, i.e.,

$$\psi(u) \not< \psi(v) \Leftrightarrow \exists i \in I: \psi_i(u) \geq \psi_i(v).$$

It was also proved in [1] that for the two-stage vector-vector choice procedure an equivalent pair-dominant choice procedure exists if and only if the alternatives from A do not make up $\varphi$, $\psi$-uncoordinated triads of the first type.

As for a pair-dominant choice procedure, it is necessary to define its levels. It was defined in [1] that a pair-dominant choice procedure has the level $\hat{1}$ if a preference $P$ on alternatives in (1) is acyclic, level $\hat{2}$ if $P$ is acyclic and transitive, i.e., it is a partial order, level $\hat{3}$ if $P$ is acyclic, transitive and negatively transitive, i.e., it is a weak order.

It turned out that for the two-stage vector-vector choice procedure to generate a choice function of level $\hat{2}$ (level $\hat{3}$), it is necessary and sufficient that no alternatives from set $A$ make up $\varphi$, $\psi$-noncoordinated triads of a second type (correspondingly $\varphi$, $\psi$-noncoordinated triads of a second and third type, and also $\psi$ , $\varphi$-noncoordinated triads of third type).

Finally, it was proved in [1] that two-stage vector-scalar choice procedure (type 3) is equivalent to

a) a pair-dominant level $\hat{1}$ choice procedure;
b) a one-stage multi-criteria extremizational choice procedure; and
c) a one-stage single-criterion extremizational choice procedure.

if and only if there are no $\varphi$-triads $\langle u,v,w \rangle$ ($\varphi$ – vector criterion) in $A$ satysfying, respectively, the following conditions a') $\psi(u) > \psi(w) \geq \psi(v)$, where $\psi$ is a scalar criterion;
b') $\psi(u) \geq \psi(w) \geq \psi(v)$ with at least one strict inequality;
c') $\psi(u) \geq \psi(w) \geq \psi(v)$.

In [1] an example of a two-stage choice procedure was provided. This two-stage choice procedure operates as follows at the first stage a Pareto set is chosen from initial set $X$ from which, at the second stage, an alternative is isolated by Euclidean metric to some «ideal» point $\omega_0 = (\omega_1, \dots, \omega_n)$ in criterial space, where $\omega_1, \dots, \omega_n$ – maximum values for each criterion. It was shown that this choice procedure can be reduced to its second stage iff

$$\forall i \in I\ \omega_i \geq \max_{x \in A} \varphi_i(x).$$

In [12] a two-stage sequential choice procedure was studied, the first stage being defined by q-Pareto multicriterial choice rule, and the second stage being defined by scalar extremization choice procedure.

Let us define a q-Pareto choice rule that was studied in [2,3]. The main idea of this rule is to choose alternatives which are dominated by no more that $q$ alternatives. Hence, it allows to avoid situations when the measurement by criterion has not been made very accurately and, consequently, it allows not to miss almost optimal elements from the chosen set of alternatives. Below we consider q-Pareto choice rule more precisely, but first we need an additional notion of upper contour set.

The upper contour set $D_i(x)$ for an alternative $x$ in the $i$-th criterion is

$$D_i(x) = \{y \epsilon A | \varphi_i(y) > \varphi_i(x)\},$$

where φ is a vector ctiterion.

In other words, the upper contour set for $x$ is a set of all alternatives which have higher values that $x$ in the criterion $i$.

The upper contour set for an alternative $x$ in a set $X$ is defined as $[D_i(x) \cap X]$.

Now, we re-define the Pareto rule as

$$x \epsilon C(X) \Leftrightarrow card\ (\cap_{i\in N}[D_i(x) \cap X]) = 0.$$

Obviously, Pareto rule is a special case of q-Pareto rule when parameter $q$ is equal to 0. The generalization of q-Pareto choice rule can be defined as

$$x \epsilon C(X) \Leftrightarrow card\ (\cap_{i\in N}[D_i(x) \cap X]) \leq q.$$

Thus, the choice rule consists in choice of alternatives which are dominated by no more than $q$ alternatives.

*Example.* Consider the following situation. $A = \{a, b, c, d, e, f, g, h, k, l, m\}$.

| **A** | *a* | *b* | *c* | *d* | *e* | *f* | *g* | *h* | *k* | *l* | *m* |
|---|---|---|---|---|---|---|---|---|---|---|---|
| $\boldsymbol{\varphi_1}$ | 1 | 3 | 5 | 0 | 5 | 4 | 4 | 5 | 2 | 4 | 1 |
| $\boldsymbol{\varphi_2}$ | 5 | 3 | 0 | 4 | 1 | 2 | 5 | 4 | 4 | 4 | 3 |

The results of applying q-Pareto choice rule depending on parameter value $q$ are the following

| ***q*** | **C(A)** |
|---|---|
| 0 | *{g,h}* |
| 1 | *{a,e,g,h}* |
| 2 | *{a,c,e,g,h,l}* |
| 3 | *{a,b,c,e,f,g,h,k,l}* |
| 4 | *{a,b,c,e,f,g,h,k,l}* |
| 5 | *{a,b,c,d,e,f,g,h,k,l}* |
| 6 | *{a,b,c,d,e,f,g,h,k,l,m}* |

In [13] it was studied which rationality conditions are satisfied for q-Pareto choice rule. It was found that q-Pareto choice rule satisfies only condition **C** (see below).

In [12] q-Pareto choice rule is used on the first stage of two-stage choice procedure. However, since the choice set of the first stage usually contains too many elements, obtained set is used as a presentation for the second stage constructed by a scalar extremization choice procedure. Thus, q-Pareto-scalar choice procedure can be presented as a superposition of two choice functions $C(\cdot)=C_2(C_1(\cdot))$, where $C_1$ is defined by q-Pareto rule, $C_2$ is defined by scalar choice procedure.

q-Pareto-scalar choice procedure was firstly studied in [12]. There were found necessary and sufficient conditions when this procedure can be reduced to the choice on scalar criterion. Two-stage choice procedure can be reduced to one-stage choice procedure if $C_1(X) \subseteq C_2(X)$, i.e., when for any subset $X \subseteq A$ the remaining set after the first stage of the rule always contains alternatives which are chosen on the second stage of the rule. In addition, necessary conditions were defined in [12] under which rationality conditions are satisfied for two-stage q-Pareto-scalar choice procedure.

Before moving to the study of the properties of two-stage choice procedures, it is necessary to consider the notion of choice functions closedness with respect to superposition operator. In order to do it some definitions were given in [1].

Any particular choice function *C(·)* can be treated as a point of the abstract space $\mathscr{C}$, consisting of all the possible choice functions on *A*.

A set of choice functions $C_1(\cdot),\dots,C_n(\cdot)$ is called a functional profile and denoted by $\{C_i(\cdot)\}$.

Let a domain $Q \subseteq \mathscr{C}$ and number *n* be fixed. Any mapping from $Q^n = \underbrace{Q \times Q \times \ldots \times Q}_{n}$ times into $\mathscr{C}$ is called a functional operator. A functional operator assigns to any profile $\{C_i(\cdot)\}$ consisting of *n* choice functions from Q a single choice function denoted by $C^*(\cdot)=F(\{C_i(\cdot)\})$, where F stands for the functional operator.

As usual, the set $Q^n$ is called the domain of F, and the set of all functions $C^*(\cdot)\epsilon\mathscr{C}$ such that $\exists\{C_i(\cdot)\} \in Q^n: F(\{C_i(\cdot)\}) = C^*(\cdot)$ as the range of operator.

The set $Q_{def} \equiv Q$ is called the domain of the functional operator $F: Q^n \to \mathscr{C}$. Any domain $Q' \subseteq \mathscr{C}$ such that $\forall\{C_i(\cdot)\}\epsilon Q^n\ F(\{C_i(\cdot)\})\epsilon Q'$ will be called the range of F($Q' = Q_{def}$). According to the definition, the domain $Q_{def}$ of the functional operator $F$ is any subset from $\mathscr{C}$ which includes the range *F*.

If $C^*(\cdot) = F(\{C_i(\cdot)\})\epsilon Q$ for any functional profile $\{C_i(\cdot)\}\epsilon Q^n$, the domain Q is called closed with respect to the operator *F*.

Now introduce functional operator of "superposition" ($F^{\odot}$):

$$F^{\odot}: \forall X \in 2^A\ C^*(X) = C_2\big(C_1(X)\big).$$

Let us define in which cases superposition, as applied to the choice functions $C_1(\cdot)$ and $C_2(\cdot)$ from different domains of $\mathscr{C}$, 'preserves' the result, i.e., the function $C^*(\cdot)$ in any of the domains. With such an approach, it is only natural to speak about 'superposition' mapping the domains $Q_1$ and $Q_2$ to which, respectively, $C_1(\cdot)$ and $C_2(\cdot)$ belong into the domain $Q'$, denoting it symbolically as $\mathfrak{R}^{\odot}(Q_1,Q_2)= Q'$.

By closedness (conditional and unconditional) of domains, we consider the following. The domain $Q_1 \subseteq \mathscr{C}$ is called *conditionally closed* at the first superposition stage at fixed domain $Q_2$ of the second stage if $\mathfrak{R}^{\odot}(Q_1, Q_2) \subseteq Q_1$. Similarly, the domain $Q_2 \subseteq \mathscr{C}$ is called conditionally closed at the second superposition stage at fixed domain $Q_1$ of the first stage if $\mathfrak{R}^{\odot}(Q_1, Q_2) \subseteq Q_2$. Finally, the domain $Q$ is called unconditionally closed with respect to the 'superposition' operator if $\mathfrak{R}^{\odot}(Q_1, Q) \subseteq Q'$.

Theorem 1 [1]. There are the following conditions of domains closedness in $\mathscr{C}$ with respect to the superposition operator ⊙:

1°. The domain **ACA** is closed (unconditionally) with respect to superposition operator, i.e., $\mathfrak{R}^{\odot}(K, K) \subseteq \boldsymbol{ACA}$; none of the domains Q = **H, C, O,** $\boldsymbol{H \cap C}$**,** $\boldsymbol{H \cap O}$**,** $\boldsymbol{C \cap O}$**,** $\boldsymbol{H \cap C \cap O}$ us closed unconditionally with respect to the superposition operator;

2°. For $Q_1$=**ACA**, the domains $Q_2$ = ***H, C, H***∩**C**, $\widehat{\boldsymbol{O}}, \widehat{\boldsymbol{H}} \cap \widehat{\boldsymbol{O}}, \widehat{\boldsymbol{C}} \cap \widehat{\boldsymbol{O}}$ and $\widehat{\boldsymbol{H}} \cap \widehat{\boldsymbol{C}} \cap \widehat{\boldsymbol{O}}$ are conditionally closed at the second superposition stage. For $Q_1$=**ACA**, the domains $Q_2$ = ***O, H***∩**O**, ***C***∩***O***, ***H***∩***C***∩***O*** are not closed. None of the domains $Q_2$=***H*** or ***C*** or... or ***H***∩***C***∩***O*** is closed for none $Q_1$ = ***H***, or ***C***, or... or ***H***∩***C***∩***O***, where the notion of $\widehat{\mathrm{H}}, \widehat{\mathrm{C}}, \widehat{\mathrm{O}}$ is used for a class of choice functions that satisfies current condition (**H**,**C**,**O**) as well as the condition of choice non-emptyness.

3°. None of the fundamental domains and their intersections **H,C, …,** ***H***∩***C***∩***O*** is conditionally closed as $Q_1$, at the first superposition stage for $Q_2$ = **ACA** (and all the more so for $Q_2$ = ***H***, **C, ...,** ***H***∩***C***∩***O*** at the second one).

*Remark*. Although generally going out of the classical domain **H**∩**C**, the superposition of two classical choice functions may be shown to stay within the fundamental domain **C**. It means that $\mathfrak{R}^{\odot}(\mathbf{H} \cap \mathbf{C}, \mathbf{H} \cap \mathbf{C}) \subset \mathbf{C}$.

Thus, the results on choice functions closedness with respect to superposition operator solve the inverse problem that helps to determine will the function resulting from application of this operator to the elements of the set satisfy the rationality conditions ***H***, ***C*** or ***O***.

Now let us study the properties of two-stage choice procedures.

## Two-stage superposition choice procedures

We consider the two-stage superposition choice procedures based on scoring rules, rules, using majority relation, value function and tournament matrix. A full description of studied choice procedures is given in [11].

A list of two-stage superposition choice procedures is provided in Table 1.

**Table 1. Two-stage choice procedures**

| № | Stage 1 | Stage 2 |
|---|---|---|
| **1-121** | Scoring rules (11 procedures) | Scoring rules (11 procedures) |
| **122-231** | | Rules, using majority relation (10 procedures) |
| **232-286** | | Rules, using value function (5 procedures) |
| **287-308** | | Rules, using tournament matrix (2 procedures) |
| **309-418** | Rules, using majority relation (10 procedures) | Scoring rules (11 procedures) |
| **419-518** | | Rules, using majority relation (10 procedures) |
| **519-568** | | Rules, using value function (5 procedures) |
| **569-588** | | Rules, using tournament matrix (2 procedures) |
| **589-643** | Rules, using value function (5 procedures) | Scoring rules (11 procedures) |
| **644-693** | | Rules, using majority relation (10 procedures) |
| **694-718** | | Rules, using value function (5 procedures) |
| **719-728** | | Rules, using tournament matrix (2 procedures) |
| **729-750** | Rules, using tournament matrix (2 procedures) | Scoring rules (11 procedures) |
| **751-770** | | Rules, using majority relation (10 procedures) |
| **771-780** | | Rules, using value function (5 procedures) |
| **781-784** | | Rules, using tournament matrix (2 procedures) |

Thus, 784 two-stage procedures of 16 different types are studied.

As properties of two-stage choice procedures should be studied separately, let us assign an identification number to each two-stage choice procedure by the following formula

$$id = 28 * (i - 1) + j,$$

where *id* is an identification number of two-stage choice procedure, *i* is a number of choice procedure from the list above used on the first stage, *j* is a number of choice procedure from the list above used on the second stage.

Before proceeding to the study of the properties of two-stage choice procedures, it is necessary to make some notes.

*Note 1*. Table 2 provides a list of two-stage procedures that does not make any sense, i.e., the second stage of which does not change the choice.

**Table 2. Two-stage choice procedures that does not make any sense («…» is any studied choice procedure)**

| Two-stage procedure | | |
|---|---|---|
| **№** | **Stage 1** | **Stage 2** |
| 1-28 | Simple majority rule | … |
| 113-140 | Run-off procedure | … |
| 141-168 | Hare rule (Ware procedure) | … |
| 231 | Inverse Borda rule | Borda rule |
| 233 | Inverse Borda rule | Inverse Borda rule |
| 234 | Inverse Borda rule | Nanson rule |
| 259 | Nanson rule | Borda rule |
| 261 | Nanson rule | Inverse Borda rule |
| 262 | Nanson rule | Nanson rule |
| 281-308 | Coombs procedure | … |
| 320 | Minimal dominant set | Minimal dominant set |
| 348 | Minimal undominant set | Minimal dominant set |
| 349 | Minimal undominant set | Minimal undominant set |
| 505-532 | Condorcet winner | … |
| 533 | Core | Simple majority rule |
| 539 | Core | Borda rule |
| 540 | Core | Black procedure |
| 541 | Core | Inverse Borda rule |
| 542 | Core | Nanson rule |
| 544 | Core | Minimal dominant set |
| 545 | Core | Minimal undominant set |
| 546 | Core | Minimal weakly stable set |
| 547 | Core | Fishburn rule |
| 548 | Core | Uncovered set I |
| 549 | Core | Uncovered set II |
| 550 | Core | Richelson rule |
| 552 | Core | Core |
| 553 | Core | k-stable set (k>1) |
| 555 | Core | Copeland rule 1 |
| 556 | Core | Copeland rule 2 |
| 557 | Core | Copeland rule 3 |
| 559 | Core | Minimax procedure |
| 560 | Core | Simpson procedure |
| **Total number of choice procedures - 168** | | |

The simple majority rule, the run-off procedure, the Hare procedure, the Coombs procedure and the Condorcet winner choose a single best alternative, consequently, two-stage choice procedures 1-28, 113-168, 281-308, 505-532 which use such procedures on the first stage do not make any sense as they do not affect the final choice.

Also two-stage procedures 231, 233, 234, 259, 261, 262 do not make any sense as the first stage of such procedures gives a set of alternatives with the same Borda count. Similarly, two-stage procedures 283, 291 do not make any sense as the first stage of such procedures gives a set of alternatives which are considered as the worst by the same number of criteria.

Two-stage procedures 320, 348, 349 do not make any sense as the second stage of such procedures does not change the choice. As for two-stage choice procedure 553, it does not make any sense, as the result of such procedure is always an empty set.

Finally, two-stage procedures 533, 539-542, 544-550, 552, 555-557,559-560 do not make any sense as $\forall x, y \in C_1(\vec{P}_X, X) \rightarrow card\{i \in N | xP_i y\} = card\{i \in N | yP_i x\}$, where $C_1(\vec{P}_X, X)$ is a first-stage procedure (the core).

Thus, choice procedures provided in Table 2 are excluded from further consideration.

*Note 2.* Properties of two-stage choice procedures which use the Black procedure on the first stage completely coincide with properties of two-stage choice procedures which use the Borda rule if there is no Condorcet winner.

*Note 3.* Table 3 provides a list of two-stage procedures which are equivalent to existing choice procedures.

**Table 3. Two-stage procedures equivalent to existing choice procedure**

| Two-stage choice procedure | | | What procedure is equivalent to |
|---|---|---|---|
| **№** | **Stage 1** | **Stage 2** | |
| 309 | Minimal dominant set | Simple majority rule | Condorcet winner |
| 321 | Minimal dominant set | Minimal undominant set | Minimal undominant set |
| 322 | Minimal dominant set | Minimal weakly stable set | Minimal weakly stable set |
| 323 | Minimal dominant set | Fishburn rule | Fishburn rule |
| 324 | Minimal dominant set | Uncovered set I | Uncovered set I |
| 325 | Minimal dominant set | Uncovered set II | Uncovered set II |
| 326 | Minimal dominant set | Richelson rule | Richelson rule |
| 327 | Minimal dominant set | Condorcet winner | Condorcet winner |
| 328 | Minimal dominant set | Core | Core |
| 331 | Minimal dominant set | Copeland rule 1 | Copeland rule 1 |
| 332 | Minimal dominant set | Copeland rule 2 | Copeland rule 2 |
| 333 | Minimal dominant set | Copeland rule 3 | Copeland rule 3 |
| 337 | Minimal undominant set | Simple majority rule | Core (single chosen alternative) |
| 350 | Minimal undominant set | Minimal weakly stable set | Minimal weakly stable set |
| 355 | Minimal undominant set | Condorcet winner | Condorcet winner |
| 356 | Minimal undominant set | Core | Core |
| 393 | Fishburn rule | Simple majority rule | Core (single chosen alternative) |
| 411 | Fishburn rule | Condorcet winner | Core (single chosen alternative) |
| 421 | Uncovered set I | Simple majority rule | Core (single chosen alternative) |
| 439 | Uncovered set I | Condorcet winner | Core (single chosen alternative) |
| 449 | Uncovered set II | Simple majority rule | Core (single chosen alternative) |
| 467 | Uncovered set II | Condorcet winner | Core (single chosen alternative) |
| 477 | Richelson rule | Simple majority rule | Core (single chosen alternative) |
| 495 | Richelson rule | Condorcet winner | Core (single chosen alternative) |
| 551 | Core | Condorcet winner | Condorcet winner |
| **Total number of choice procedures - 25** | | | |

Two-stage choice procedures provided in previous table are equivalent to some existing choice procedures. Thus, their properties fully coincided with properties of such existing procedures. However, these two-stage procedures are not excluded from further consideration as the computational complexity of some of them can be lower than the complexity of existing procedures.

Thus, it remains to study properties of 591 two-stage choice procedures.

## A study of the properties of two-stage choice procedures

A list of studied normative properties is given in [11].

A study of the properties is conducted as follows. If a two-stage choice procedure does not satisfy given normative condition, a counter-example is provided. On the country, if a two-stage choice procedure satisfies given normative condition a necessary proof is followed. The study of the properties of two-stage choice procedures is based on the study of the properties of multi-criteria choice procedures which is done in [11].

The results of the study of the properties of 591 two-stage choice procedures are given in Theorem 2.

**Theorem 2.** Information on which choice procedures satisfy given normative conditions is provided in Table 4.

**Table 4. Properties of two-stage choice procedures («+» - choice procedure satisfies given normative condition, «-» - choice procedure does not satisfy given normative condition)**

| Two-stage choice procedure | | | | Normative conditions | | | | | | | |
|---|---|---|---|---|---|---|---|---|---|---|---|
| | | | | Rationality conditions | | | | Monotonicity conditions | | | |
| Stage 1 | | Stage 2 | | Heredity condition (H) | Concordance condition (C) | Outcast condition (O) | Heredity condition (H) | Monotonicity condition 1 | Monotonicity condition 2 | Strict monotonicity condition | Non-compensatory condition |
| **№** | **Name** | **№** | **Name** | | | | | | | | |
| **2**<br>**3**<br>**4**<br>**12**<br>**13**<br>**15**<br>**16**<br>**17**<br>**18**<br>**20**<br>**23**<br>**24**<br>**25**<br>**26**<br>**27**<br>**28** | Plurality rule<br>Inverse plurality rule<br>q-Approval rule (q>1)<br>Minimal dominant set<br>Minimal undominant set<br>Fishburn rule<br>Uncovered set I<br>Uncovered set II<br>Richelson rule<br>Core<br>Copeland rule 1<br>Copeland rule 2<br>Copeland rule 3<br>Super-threshold rule (threshold depends on X)<br>Minimax procedure<br>Simpson procedure | **5**<br>**6**<br>**11** | Run-off procedure<br>Hare rule (Ware procedure)<br>Coombs procedure | - | - | - | - | - | + | - | - |
| **9**<br>**10**<br>**11**<br>**14**<br>**21** | Inverse Borda rule<br>Nanson rule<br>Coombs procedure<br>Minimal weakly stable set<br>k-stable set (k>1) | **1**<br>**5**<br>**6**<br>**11**<br>**19** | Simple majority rule<br>Run-off procedure<br>Hare rule (Ware procedure)<br>Coombs procedure<br>Condorcet winner | - | - | - | - | - | + | - | - |
| **7**<br>**8**<br>**22** | Borda rule<br>Black procedure<br>Threshold rule | **1**<br>**5**<br>**6** | Simple majority rule<br>Run-off procedure<br>Hare rule (Ware procedure) | - | - | - | - | + | + | - | - |

| Two-stage choice procedure | | | | Normative conditions | | | | | | | |
|---|---|---|---|---|---|---|---|---|---|---|---|
| Stage 1 | | Stage 2 | | Rationality conditions | | | | Monotonicity conditions | | | Non-compensatory condition |
| № | Name | № | Name | Heredity condition (H) | Concordance condition (C) | Outcast condition (O) | Heredity condition (H) | Monotonicity condition 1 | Monotonicity condition 2 | Strict monotonicity condition | |
| | | **19** | Condorcet winner | | | | | | | | |
| **2**<br>**3**<br>**4**<br>**23**<br>**24**<br>**25**<br>**26**<br>**27**<br>**28** | Plurality rule<br>Inverse plurality rule<br>q-Approval rule (q>1)<br>Copeland rule 1<br>Copeland rule 2<br>Copeland rule 3<br>Super-threshold rule (threshold depends on X)<br>Minimax procedure<br>Simpson procedure | **1**<br>**19** | Simple majority rule<br>Condorcet winner | - | - | - | - | + | + | - | - |
| **13**<br>**20** | Minimal undominant set<br>Core | **1** | Simple majority rule | + | + | - | - | + | + | - | - |
| **15**<br>**16**<br>**17**<br>**18** | Fishburn rule<br>Uncovered set I<br>Uncovered set II<br>Richelson rule | **1**<br>**19** | Simple majority rule<br>Condorcet winner | + | + | - | - | + | + | - | - |
| **7**<br>**8**<br>**22**<br>**23**<br>**24**<br>**25** | Borda rule<br>Black procedure<br>Threshold rule<br>Copeland rule 1<br>Copeland rule 2<br>Copeland rule 3 | **9**<br>**10**<br>**21** | Inverse Borda rule<br>Nanson rule<br>k-stable set (k>1) | - | - | - | - | + | - | - | - |
| **13** | Minimal undominant set | **15**<br>**17**<br>**21**<br>**23**<br>**24**<br>**25** | Fishburn rule<br>Uncovered set II<br>k-stable set (k>1)<br>Copeland rule 1<br>Copeland rule 2<br>Copeland rule 3 | - | - | - | - | + | - | - | - |
| **15** | Fishburn rule | **12** | Minimal dominant set | - | - | - | - | + | - | - | - |

| Two-stage choice procedure | | | | Normative conditions | | | | | | | |
|---|---|---|---|---|---|---|---|---|---|---|---|
| | | | | Rationality conditions | | | | Monotonicity conditions | | | Non-compensatory condition |
| Stage 1 | | Stage 2 | | Heredity condition (H) | Concordance condition (C) | Outcast condition (O) | Heredity condition (H) | Monotonicity condition 1 | Monotonicity condition 2 | Strict monotonicity condition | |
| **№** | **Name** | **№** | **Name** | | | | | | | | |
| **16**<br>**18** | Richelson rule<br>Uncovered set I | **16**<br>**18** | Uncovered set I<br>Richelson rule | | | | | | | | |
| **26**<br>**27**<br>**28** | Super-threshold rule (threshold depends on X)<br>Minimax procedure<br>Simpson procedure | **12**<br>**13**<br>**20** | Minimal dominant set<br>Minimal undominant set<br>Core | - | - | - | - | + | - | - | - |
| **2**<br>**3**<br>**4**<br>**7**<br>**8**<br>**22**<br>**23**<br>**24**<br>**25** | Plurality rule<br>Inverse plurality rule<br>q-Approval rule (q>1)<br>Borda rule<br>Black procedure<br>Threshold rule<br>Copeland rule 1<br>Copeland rule 2<br>Copeland rule 3 | **2**<br>**3**<br>**4**<br>**7**<br>**8**<br>**12**<br>**13**<br>**14**<br>**15**<br>**16**<br>**17**<br>**18**<br>**20**<br>**22**<br>**23**<br>**24**<br>**25**<br>**26**<br>**27**<br>**28** | Plurality rule<br>Inverse plurality rule<br>q-Approval rule (q>1)<br>Borda rule<br>Black procedure<br>Minimal dominant set<br>Minimal undominant set<br>Minimal weakly stable set<br>Fishburn rule<br>Uncovered set I<br>Uncovered set II<br>Richelson rule<br>Core<br>Threshold rule<br>Copeland rule 1<br>Copeland rule 2<br>Copeland rule 3<br>Super-threshold rule (threshold depends on X)<br>Minimax procedure<br>Simpson procedure | - | - | - | - | + | - | - | - |
| **2**<br>**3**<br>**4** | Plurality rule<br>Inverse plurality rule<br>q-Approval rule (q>1) | **21** | k-stable set (k>1) | - | - | - | - | + | - | - | - |
| **13** | Minimal undominant set | **16** | Uncovered set I | - | - | - | - | + | - | - | - |

| Two-stage choice procedure | | | | Normative conditions | | | | | | | |
|---|---|---|---|---|---|---|---|---|---|---|---|
| | | | | Rationality conditions | | | | Monotonicity conditions | | | |
| | | | | Heredity condition (H) | Concordance condition (C) | Outcast condition (O) | Heredity condition (H) | Monotonicity condition 1 | Monotonicity condition 2 | Strict monotonicity condition | Non-compensatory condition |
| **Stage 1** | | **Stage 2** | | | | | | | | | |
| **№** | **Name** | **№** | **Name** | | | | | | | | |
| **17** | Uncovered set II | **18** | Richelson rule | | | | | | | | |
| **16** | Uncovered set I | **20** | Core | - | - | - | - | + | - | - | - |
| **12**<br>**13** | Minimal dominant set<br>Minimal undominant set | **26**<br>**27**<br>**28** | Super-threshold rule (threshold depends on X)<br>Minimax procedure<br>Simpson procedure | - | - | - | - | + | - | - | - |
| **20** | Core | **7**<br>**8**<br>**9**<br>**10**<br>**12**<br>**13**<br>**14**<br>**15**<br>**16**<br>**17**<br>**18**<br>**23**<br>**24**<br>**25**<br>**27**<br>**28** | Borda rule<br>Black procedure<br>Inverse Borda rule<br>Nanson rule<br>Minimal dominant set<br>Minimal undominant set<br>Minimal weakly stable set<br>Fishburn rule<br>Uncovered set I<br>Uncovered set II<br>Richelson rule<br>Copeland rule 1<br>Copeland rule 2<br>Copeland rule 3<br>Minimax procedure<br>Simpson procedure | + | + | - | - | + | + | - | - |
| **17** | Uncovered set II | **12** | Minimal dominant set | - | + | - | - | + | - | - | - |
| **9**<br>**10**<br>**14**<br>**21** | Inverse Borda rule<br>Nanson rule<br>Minimal weakly stable set<br>k-stable set (k>1) | **12**<br>**13**<br>**16**<br>**18**<br>**20** | Minimal dominant set<br>Minimal undominant set<br>Uncovered set I<br>Richelson rule<br>Core | - | - | - | - | - | - | - | - |
| **15**<br>**17** | Fishburn rule<br>Uncovered set II | **13**<br>**20** | Minimal undominant set<br>Core | - | - | - | - | - | - | - | - |

| Two-stage choice procedure | | | | Normative conditions | | | | | | | |
|---|---|---|---|---|---|---|---|---|---|---|---|
| | | | | Rationality conditions | | | | Monotonicity conditions | | | Non-compensatory condition |
| Stage 1 | | Stage 2 | | Heredity condition (H) | Concordance condition (C) | Outcast condition (O) | Heredity condition (H) | Monotonicity condition 1 | Monotonicity condition 2 | Strict monotonicity condition | |
| **№** | **Name** | **№** | **Name** | | | | | | | | |
| **18** | Richelson rule | | | | | | | | | | |
| **9**<br>**10**<br>**14**<br>**15**<br>**16**<br>**17**<br>**18**<br>**21**<br>**26**<br>**27**<br>**28** | Inverse Borda rule<br>Nanson rule<br>Minimal weakly stable set<br>Fishburn rule<br>Uncovered set I<br>Uncovered set II<br>Richelson rule<br>k-stable set (k>1)<br>Super-threshold rule (threshold depends on X)<br>Minimax procedure<br>Simpson procedure | **2**<br>**4**<br>**8**<br>**14**<br>**15**<br>**17**<br>**21**<br>**22**<br>**23**<br>**24**<br>**25**<br>**26**<br>**27**<br>**28** | Plurality rule<br>q-Approval rule (q>1)<br>Black procedure<br>Minimal weakly stable set<br>Fishburn rule<br>Uncovered set II<br>k-stable set (k>1)<br>Threshold rule<br>Copeland rule 1<br>Copeland rule 2<br>Copeland rule 3<br>Super-threshold rule (threshold depends on X)<br>Minimax procedure<br>Simpson procedure | - | - | - | - | - | - | - | - |
| **12**<br>**13** | Minimal dominant set<br>Minimal undominant set | **2**<br>**3**<br>**4**<br>**7**<br>**8**<br>**22** | Plurality rule<br>Inverse plurality rule<br>q-Approval rule (q>1)<br>Borda rule<br>Black procedure<br>Threshold rule | - | - | - | - | - | - | - | - |
| **14**<br>**15**<br>**16**<br>**17**<br>**18**<br>**21**<br>**26**<br>**27** | Minimal weakly stable set<br>Fishburn rule<br>Uncovered set I<br>Uncovered set II<br>Richelson rule<br>k-stable set (k>1)<br>Super-threshold rule (threshold depends on X)<br>Minimax procedure | **3**<br>**7**<br>**9**<br>**10** | Inverse plurality rule<br>Borda rule<br>Inverse Borda rule<br>Nanson rule | - | - | - | - | - | - | - | - |

| Two-stage choice procedure | | | | Normative conditions | | | | | | | |
|---|---|---|---|---|---|---|---|---|---|---|---|
| | | | | Rationality conditions | | | | Monotonicity conditions | | | |
| Stage 1 | | Stage 2 | | Heredity condition (H) | Concordance condition (C) | Outcast condition (O) | Heredity condition (H) | Monotonicity condition 1 | Monotonicity condition 2 | Strict monotonicity condition | Non-compensatory condition |
| **№** | **Name** | **№** | **Name** | | | | | | | | |
| **28** | Simpson procedure | | | | | | | | | | |
| **9**<br>**10** | Inverse Borda rule<br>Nanson rule | **3**<br>**11** | Inverse plurality rule | - | - | - | - | - | - | - | - |
| **16** | Uncovered set I | **13** | Minimal undominant set | - | - | - | - | - | - | - | - |
| **20** | Core | **2**<br>**3**<br>**4**<br>**22**<br>**26** | Plurality rule<br>Inverse plurality rule<br>q-Approval rule (q>1)<br>Threshold rule<br>Super-threshold rule (threshold depends on X) | - | - | - | - | - | - | - | - |
| **2**<br>**3**<br>**4**<br>**12**<br>**13**<br>**23**<br>**24**<br>**25** | Plurality rule<br>Inverse plurality rule<br>q-Approval rule (q>1)<br>Minimal dominant set<br>Minimal undominant set<br>Copeland rule 1<br>Copeland rule 2<br>Copeland rule 3 | **9**<br>**10** | Inverse Borda rule<br>Nanson rule | - | - | - | - | - | - | - | - |
| **26**<br>**27**<br>**28** | Super-threshold rule (threshold depends on X)<br>Minimax procedure<br>Simpson procedure | **16**<br>**18** | Uncovered set I<br>Richelson rule | - | - | - | - | - | - | - | - |

The proof of the theorem is provided in Appendix 1.

## Computational complexity of choice procedures

A computational complexity of choice procedures used in two-stage superpositions is provided in Table 5.

**Table 5. Theoretical computational complexity of existing choice procedures (M – cardinality of initial set of alternatives, n – number of criteria, *q,k* – parameters of choice procedures, *l,d* – parameters which depend on initial set of alternatives, $d \geq 1, 1 \leq l \leq M$)**

| № | Name of choice procedure | Computational complexity (theoretical) | Maximum number of remaining alternatives after applying the choice procedure |
|---|---|---|---|
| **1** | Simple majority rule | $O(M \cdot n)$ | 1 |
| **2** | Plurality rule | $O(M \cdot n)$ | $n$ |
| **3** | Inverse plurality rule | $O(M \cdot n)$ | $M - n$ |
| **4** | q-Approval rule (q>1) | $O(M \cdot n \cdot \log_2(q+1))$ | $q \cdot n$ |
| **5** | Run-off procedure | $O(M \cdot n)$ | 1 |
| **6** | Hare rule (Ware procedure) | $O(M \cdot n)$ | 1 |
| **7** | Borda rule | $O(M \cdot \log_2(M) \cdot n)$ | $M$ |
| **8** | Black procedure | $O(M \cdot \log_2(M) \cdot n)$ | $M$ |
| **9** | Inverse Borda rule | $O(M^2 \cdot n)$ | $M$ |
| **10** | Nanson rule | $O(M^2 \cdot n)$ | $M$ |
| **11** | Coombs procedure | $O(M^2 \cdot n)$ | $n$ |
| **12** | Minimal dominant set | $O(M^{2,37 \cdot d})$ | $M$ |
| **13** | Minimal undominant set | $O(M^{2,37 \cdot d})$ | $M$ |
| **14** | Minimal weakly stable set | $O(M^2 \cdot n + \sum_{i=1}^{l} C_M^i \cdot (M - i))$ | $M$ |
| **15** | Fishburn rule | $O(M^3)$ | $M$ |
| **16** | Uncovered set I | $O(M^3)$ | $M$ |
| **17** | Uncovered set II | $O(M^3)$ | $M$ |
| **18** | Richelson rule | $O(M^3)$ | $M$ |
| **19** | Condorcet winner | $O(M \cdot n)$ | 1 |
| **20** | Core | $O(M^2 \cdot n)$ | $M$ |
| **21** | k-stable set (k>1) | $O(M^{2,37 \cdot k})$ | $M$ |
| **22** | Threshold rule | $O(M \cdot n)$ | $M$ |
| **23** | Copeland rule 1 | $O(M^2 \cdot n)$ | $M$ |
| **24** | Copeland rule 2 | $O(M^2 \cdot n)$ | $M$ |
| **25** | Copeland rule 3 | $O(M^2 \cdot n)$ | $M$ |
| **26** | Super-threshold rule | $O(M)$ | $M$ |
| **27** | Minimax procedure | $O(M^2 \cdot n)$ | $M$ |
| **28** | Simpson procedure | $O(M^2 \cdot n)$ | $M$ |

A computational complexity of choice procedures was calculated by the author of the paper.

Based on information provided in Table 5 we can divide all two-stage procedures in several groups in accordance with their computational complexity. The results are provided in Table 6.

**Table 6. A computational complexity of two-stage choice procedures («…» is any choice procedure)**

| Two-stage choice procedure | |
|---|---|
| **Stage 1** | **Stage 2** |
| **Choice procedures with a low computational complexity** | |
| Plurality rule<br>Threshold rule<br>q-Approval rule (q>1) | **…** |
| Inverse plurality rule<br>Super-threshold rule<br>Borda rule<br>Black procedure | Simple majority rule<br>Run-off procedure<br>Hare rule (Ware procedure)<br>Borda rule<br>Black procedure<br>Condorcet winner<br>Plurality rule<br>Threshold rule<br>Inverse plurality rule<br>q-Approval rule (q>1)<br>Super-threshold rule |
| **Computational complexity depends on initial set of alternatives** | |
| Inverse plurality rule<br>Super-threshold rule<br>Borda rule<br>Black procedure | Inverse Borda rule<br>Nanson rule<br>Core<br>Copeland rules 1-3<br>Minimax procedure<br>Simpson procedure<br>Coombs procedure |
| **Choice procedures with average computational complexity** | |
| Inverse Borda rule<br>Nanson rule<br>Core<br>Copeland rules 1-3<br>Minimax procedure<br>Simpson procedure | Simple majority rule<br>Run-off procedure<br>Hare rule (Ware procedure)<br>Borda rule<br>Black procedure<br>Condorcet winner<br>Plurality rule<br>Threshold rule<br>Inverse plurality rule<br>q-Approval rule (q>1)<br>Super-threshold rule<br>Inverse Borda rule<br>Nanson rule<br>Core<br>Copeland rules 1-3<br>Minimax procedure<br>Simpson procedure<br>Coombs procedure |
| Coombs procedure | **…** |
| **Choice procedures with a high computational complexity** | |
| Inverse plurality rule<br>Super-threshold rule<br>Borda rule<br>Black procedure<br>Inverse Borda rule | Minimal dominant set<br>Minimal undominant set<br>Minimal weakly stable set<br>Fishburn rule<br>Uncovered set I, II |

| Two-stage choice procedure | |
|---|---|
| **Stage 1** | **Stage 2** |
| Nanson rule<br>Core<br>Copeland rules 1-3<br>Minimax procedure<br>Simpson procedure | Richelson rule<br>k-stable set (k>1) |
| Minimal dominant set<br>Minimal undominant set<br>Minimal weakly stable set<br>Fishburn rule<br>Uncovered set I, II<br>Richelson rule<br>k-stable set (k>1) | … |

## The run-time complexity for two-stage superposition choice procedures

Let the initial set $M_1$ contains 300 thousands of alternatives and the total number of criteria is equal to 10. Suppose, a computer can process around 3 billion instructions per second. Consider situations when the remaining set $M_2$ after applying the first-stage choice procedure contains 10, 50 and 100 thousands of alternatives.

The run-time complexity of two-stage choice procedures from different groups (according to Table 6) is provided in Table 7.

**Table 7. The run-time complexity for two-stage choice procedures ($M_2$ – number of alternatives remained after applying the first stage of the choice procedure)**

| Two-stage choice procedures | | Run time | | |
|---|---|---|---|---|
| **Stage 1** | **Stage 2** | | | |
| **Choice procedures with a low computational complexity** | | **$M_2$=10000** | **$M_2$=50000** | **$M_2$=100000** |
| Plurality rule | Uncovered set I | ≈50 ms | | |
| Inverse plurality rule | Condorcet winner | ≈52 ms | ≈58 ms | ≈66 ms |
| Super-threshold rule | Threshold rule | ≈6 ms | ≈13 ms | ≈21 ms |
| Borda rule | Plurality rule | ≈910 ms | ≈918 ms | ≈926 ms |
| **Computational complexity depends on initial set of alternatives** | | | | |
| Inverse plurality rule | Minimax procedure | ≈16 s | ≈7 min | ≈27 min |
| **Choice procedures with average computational complexity** | | | | |
| Inverse Borda rule | Simple majority rule | ≈4 h 12 min | | |
| Minimax procedure | Simpson procedure | ≈4 h 15 min | ≈4 h 21 min | ≈4 h 36 min |
| **Choice procedures with a high computational complexity** | | | | |
| Inverse plurality rule | Richelson rule | ≈46 h 20 min | ≈241 days | ≈5 years |
| Uncovered set I | Borda rule | ≈140 years | | |

The results obtained from Tables 6 and 7 give us information on which procedures can be applied when we deal with Big Data and which procedures cannot be applied in such problems as they are not allow to obtain any results in a sufficient time.

# Conclusion

We have studied the properties of 592 two-stage choice procedures, which can be used in various multi-criteria problems. It was defined which choice procedures satisfy given normative conditions, showing how a final choice is changed due to the changes of preferences or a set of feasible alternatives. Such information leads to a better understanding of different choice procedures and how stable and sensible is a set of alternatives obtained after applying some choice procedure.

The results show that most of the two-stage procedures do not satisfy any normative conditions. Only some of them satisfy the Monotonicity conditionon 1. More detailed information is provided in Table 4.

To compute run-time complexity of choice procedures the average computational complexity was used. All choice procedures were divided into different groups (see Table 6). It was shown that two-stage choice procedures which use choice procedures with a high computational complexity on the first stage require more time than other procedures. It means that such procedures are not recommended to use in applications to Big Data. Two-stage choice procedures which use on the first stage choice procedures with a low computational complexity and on the second stage - with a high computational complexity can be used in applications to Big Data, however, their application depends on the number of alternatives remained after the first stage. Two-stage choice procedures, which use on both stages choice procedures with a low computational complexity, can be used in applications to Big Data with no restrictions.

# Appendix 1. Properties of two-stage choice procedures

## Two-stage choice procedure 29. 'Plurality rule – Simple majority rule'

### *29.1 Heredity condition (H)*

Let $X = \{a, b, c\}$ and the profile $\vec{P}_X$ is the following

| $P_1$ | $P_2$ | $P_3$ | $P_4$ | $P_5$ |
|---|---|---|---|---|
| $a$ | $a$ | $c$ | $b$ | $b$ |
| $c$ | $c$ | $b$ | $a$ | $a$ |
| $b$ | $b$ | $a$ | $c$ | $c$ |

According to the two-stage choice procedure $C(\vec{P}_X, X) = \{b\}$.

Consider now the subset $X' = X\setminus\{a\}$. A contraction of a profile $\vec{P}_X$ onto a set $X'$, i.e., $\vec{P}_{X'}$, looks as

| $P_1$ | $P_2$ | $P_3$ | $P_4$ | $P_5$ |
|---|---|---|---|---|
| $c$ | $c$ | $c$ | $b$ | $b$ |
| $b$ | $b$ | $b$ | $c$ | $c$ |

According to the rule the alternative *c* will be chosen, i.e., $C\left(\vec{P}_{X'}, X'\right) = \{c\}$.

Then $C(\vec{P}_{X'}, X') \not\supseteq C(\vec{P}_X, X) \cap X'$. Thus, the condition **H** is not satisfied.

### *29.2 Concordance condition (C)*

Let $X = \{a, b, c\}$ and the profile $\vec{P}_X$ is the following

| $P_1$ | $P_2$ | $P_3$ |
|---|---|---|
| $a$ | $b$ | $c$ |
| $b$ | $a$ | $b$ |
| $c$ | $c$ | $a$ |

According to the two-stage choice procedure the choice is empty, i.e., $C(\vec{P}_X, X) = \emptyset$.

Now let us consider the subset $X' = X\setminus\{c\}$. A contraction of a profile $\vec{P}_X$ onto a set $X'$, i.e., $\vec{P}_{X'}$, looks as

| $P_1$ | $P_2$ | $P_3$ |
|---|---|---|
| $a$ | $b$ | $b$ |
| $b$ | $a$ | $a$ |

According to the two-stage choice procedure the alternative *b* will be chosen, i.e., $C\left(\vec{P}_{X'}, X'\right) = \{b\}$.

Finally, let us consider the subset $X'' = X\setminus\{a\}$. A contraction of a profile $\vec{P}_X$ onto a set $X''$, i.e., $\vec{P}_{X''}$, looks as

| $P_1$ | $P_2$ | $P_3$ |
|---|---|---|
| $b$ | $b$ | $c$ |
| $c$ | $c$ | $b$ |

According to the two-stage choice procedure the alternative *b* will be chosen, i.e., $C\left(\vec{P}_{X''}, X''\right) = \{b\}$.

Then $C\left(\vec{P}_{X'}, X'\right) \cap C\left(\vec{P}_{X''}, X''\right) = \{b\} \not\subseteq C(\vec{P}_X, X)$. Thus, the condition **C** is not satisfied.

### *29.3 Outcast condition (O)*

The condition **O** is not satisfied (see paragraph 29.1 of Appendix 1).

### *29.4 Arrow's choice axiom (ACA)*

The condition **ACA** is not satisfied since the two-stage choice procedure does not satisfy the condition **H**.

***29.5 Monotonicity condition 1***

Let $C(\vec{P}_X, X) = \{a\}$. It means that

1. $\forall x \in X\backslash\{a\}\, n^+(a, \vec{P}_X) > n^+(x, \vec{P}_X)$ or
2. $\exists y \in X\backslash\{a\}\, n^+(a, \vec{P}_X) = n^+(y, \vec{P}_X)$ and $n^+(a, \vec{P}_{X'}) > \frac{1}{2} \cdot |\vec{P}_{X'}|$, where $X'$ is a set of alternatives remained for the second stage of the choice procedure.

Consider now a profile $\vec{P'_X}$, which differs from the profile $\vec{P}_X$ only by improved position of the alternative $a$. Then

- $n^+\left(a, \vec{P'_X}\right) > n^+(a, \vec{P}_X), n^+\left(x, \vec{P'_X}\right) \leq n^+(x, \vec{P}_X)\ or$
- $n^+\left(a, \vec{P'_X}\right) = n^+(a, \vec{P}_X), n^+\left(a, \vec{P'_{X'}}\right) \geq n^+(a, \vec{P}_{X'})$ and $|\vec{P}_X| = |\vec{P'_X}|$.

Thus, $C\left(\vec{P'_X}, X\right) = \{a\}$.

Then $a \in C(\vec{P}_X, X)$ and $a \in C(\vec{P'_X}, X)$. Consequently, the Monotonicity conditionon 1 is satisfied.

***29.6 Monotonicity condition 2***

Since given two-stage choice procedure chooses no more than one best alternative, the Monotonicity conditionon 2 is not applicable to it as it considers the choice of more than two alternatives. In other words, such two-stage choice procedure obeys the Monotonicity conditionon 2 trivially.

***29.7 Strict monotonicity condition***

Let $X = \{a, b, c\}$ and the profile $\vec{P}_X$ is the following

| $P_1$ | $P_2$ | $P_3$ |
|---|---|---|
| $a$ | $a$ | $b$ |
| $c$ | $b$ | $c$ |
| $b$ | $c$ | $a$ |

According to the two-stage choice procedure the alternative *a* will be chosen, i.e., $C(\vec{P}_X, X) = \{a\}$.

Consider now a profile $\vec{P'_X}$, which differs from the profile $\vec{P}_X$ only by improved position of the alternative *c* in $P'_1$:

| $P'_1$ | $P'_2$ | $P'_3$ |
|---|---|---|
| $c$ | $a$ | $b$ |
| $a$ | $b$ | $c$ |
| $b$ | $c$ | $a$ |

According to the two-stage choice procedure the choice is empty, i.e., $C(\vec{P'_X}, X) = \emptyset$.

$$C\left(\vec{P'_X}, X\right) \neq \left[\begin{array}{c} C(\vec{P}_X, X)\ or \\ \{c\}\ or \\ C(\vec{P}_X, X) \cup \{c\}. \end{array}\right.$$

Thus, the strict monotonicity condition is not satisfied.

***29.8 Non-compensatory condition***

Let $X = \{a, b, c\}$ and the profile $\vec{P}_X$ is the following

| $P_1$ | $P_2$ | $P_3$ |
|---|---|---|
| $a$ | $a$ | $b$ |
| $b$ | $b$ | $c$ |
| $c$ | $c$ | $a$ |

According to the rule the alternative $a$ will be chosen, i.e., $C(\vec{P}_X, X) = \{a\}$.
Let us write the profile $\vec{P}_X$ in the following form

| $X$ | $\varphi_1$ | $\varphi_2$ | $\varphi_3$ |
|---|---|---|---|
| $a$ | 3 | 3 | 1 |
| $b$ | 2 | 2 | 3 |
| $c$ | 1 | 1 | 2 |

According to the non-compensatory condition the alternative $b$ is better than the alternative $a$ and the alternative $a$ is better than the alternative $c$. Thus, the non-compensatory condition is not satisfied as $\{b\} \neq C(\vec{P}_X, X)$.

## Two-stage choice procedures 30-46.

Two-stage choice procedures 30-46 do not satisfy the same conditions as choice procedures which are used on the second stage. To prove it we can use the same examples but with larger number of alternatives omitted on the first stage of the choice procedure [11].

Example. Let $X = \{a, b, c\}$ and the profile $\vec{P}_X$ is the following

| $P_1$ | $P_2$ | $P_3$ | $P_4$ | $P_5$ |
|---|---|---|---|---|
| *a* | *a* | *c* | *b* | *b* |
| *c* | *c* | *b* | *a* | *a* |
| *b* | *b* | *a* | *c* | *c* |

Let us transform this example to the following form

| $P_1$ | $P_2$ | $P_3$ | $P_4$ | $P_5$ | $P_6$ | $P_7$ | $P_8$ | $P_9$ | $P_{10}$ | $P_{11}$ |
|---|---|---|---|---|---|---|---|---|---|---|
| *d* | *e* | *f* | *g* | *h* | *a* | *b* | *c* | *a* | *b* | *c* |
| *a* | *a* | *c* | *b* | *b* | *b* | *c* | *a* | *c* | *a* | *b* |
| *c* | *c* | *b* | *a* | *a* | *c* | *a* | *b* | *b* | *c* | *a* |
| *b* | *b* | *a* | *c* | *c* | | | … | | | |
| | | | | | … | | | | | |

According to the plurality rule the alternatives *d,e,f,g,h* are omitted on the first stage of the procedure and the alternatives *a,b,c* are presented for choice on the second stage.

Thus, it is necessary to check those normative conditions of two-stage choice procedures which are satisfied for the choice procedures of the second stage.

- Two-stage choice procedures 30-32,35-36,40-56 satisfy the Monotonicity conditionon 1 (the proof follows from the paragraph 29.5 of Appendix 1 and properties of second stage choice procedures).
- Two-stage choice procedures 40,45,47,48 do not satisfy the condition **C**. To proof it the following example is used.

Let $X = \{a, b, c, d\}$ and the profile $\vec{P}_X$ is the following

| $P_1$ | $P_2$ | $P_3$ | $P_4$ | $P_5$ | $P_6$ | $P_7$ | $P_8$ | $P_9$ | $P_{10}$ | $P_{11}$ | $P_{12}$ | $P_{13}$ | $P_{14}$ |
|---|---|---|---|---|---|---|---|---|---|---|---|---|---|
| *c* | *c* | *c* | *b* | *b* | *b* | *d* | *d* | *d* | *d* | *a* | *a* | *a* | *a* |
| *a* | *a* | *a* | *d* | *d* | *d* | *b* | *b* | *c* | *c* | *b* | *b* | *b* | *c* |
| *d* | *d* | *d* | *a* | *a* | *a* | *c* | *c* | *b* | *b* | *c* | *c* | *c* | *b* |
| *b* | *b* | *b* | *c* | *c* | *c* | *a* | *a* | *a* | *a* | *d* | *d* | *d* | *d* |

Let us calculate the value $n^{+}(x, \vec{P}_X)$ for each alternative $x \in X$: $n^{+}(a, \vec{P}_X) = 4$, $n^{+}(b, \vec{P}_X) = 3$, $n^{+}(c, \vec{P}_X) = 3$, $n^{+}(d, \vec{P}_X) = 4$. According to the rule the alternatives $a$ and $d$ will be chosen, i.e., $C(\vec{P}_X, X) = \{a, d\}$.

Consider now the subset $X' = X \backslash \{d\}$. A contraction of a profile $\vec{P}_X$ onto a set $X'$, i.e., $\vec{P}_{X'}$, looks as

| $P_1$ | $P_2$ | $P_3$ | $P_4$ | $P_5$ | $P_6$ | $P_7$ | $P_8$ | $P_9$ | $P_{10}$ | $P_{11}$ | $P_{12}$ | $P_{13}$ | $P_{14}$ |
|---|---|---|---|---|---|---|---|---|---|---|---|---|---|
| $c$ | $c$ | $c$ | $b$ | $b$ | $b$ | $b$ | $b$ | $c$ | $c$ | $a$ | $a$ | $a$ | $a$ |
| $a$ | $a$ | $a$ | $a$ | $a$ | $a$ | $c$ | $c$ | $b$ | $b$ | $b$ | $b$ | $b$ | $c$ |
| $b$ | $b$ | $b$ | $c$ | $c$ | $c$ | $a$ | $a$ | $a$ | $a$ | $c$ | $c$ | $c$ | $b$ |

Let us calculate the value $n^{+}(x, \vec{P}_{X'})$ for each alternative $x \in X'$: $n^{+}(a, \vec{P}_{X'}) = 4$, $n^{+}(b, \vec{P}_{X'}) = 5$, $n^{+}(c, \vec{P}_{X'}) = 5$. Thus, the alternative $a$ will be omitted on the first stage and the alternative $b$ will be chosen on the second stage, i.e., $C(\vec{P}_{X'}, X') = \{b\}$.

Finally, let us consider the subset $X'' = X \backslash \{a\}$. A contraction of a profile $\vec{P}_X$ onto a set $X''$, i.e., $\vec{P}_{X''}$, looks as

| $P_1$ | $P_2$ | $P_3$ | $P_4$ | $P_5$ | $P_6$ | $P_7$ | $P_8$ | $P_9$ | $P_{10}$ | $P_{11}$ | $P_{12}$ | $P_{13}$ | $P_{14}$ |
|---|---|---|---|---|---|---|---|---|---|---|---|---|---|
| $c$ | $c$ | $c$ | $b$ | $b$ | $b$ | $d$ | $d$ | $d$ | $d$ | $b$ | $b$ | $b$ | $c$ |
| $d$ | $d$ | $d$ | $d$ | $d$ | $d$ | $b$ | $b$ | $c$ | $c$ | $c$ | $c$ | $c$ | $b$ |
| $b$ | $b$ | $b$ | $c$ | $c$ | $c$ | $c$ | $c$ | $b$ | $b$ | $d$ | $d$ | $d$ | $d$ |

Let us calculate the value $n^{+}(x, \vec{P}_{X''})$ for each alternative $x \in X''$: $n^{+}(b, \vec{P}_{X''}) = 6$, $n^{+}(c, \vec{P}_{X''}) = 5$, $n^{+}(d, \vec{P}_{X''}) = 4$. According to the rule the alternative $b$ will be chosen, i.e., $C(\vec{P}_{X''}, X'') = \{b\}$.

Then $C(\vec{P}_{X'}, X') \cap C(\vec{P}_{X''}, X'') = \{b\} \not\subseteq C(\vec{P}_X, X)$. Thus, the two-stage choice procedure does not satisfy the condition **C**.

- Two-stage choice procedures 40,47 do not satisfy the condition **O** (see paragraph 29.1 of Appendix 1).
- Two-stage choice procedures 47,48 do not satisfy the condition **H** (see paragraph 29.1 of Appendix 1).
- Two-stage choice procedures 48,54 do not satisfy the Monotonicity conditionon 2. To prove it the following example is used.

Let $X = \{a, b, c\}$ and the profile $\vec{P}_X$ is the following

| $P_1$ | $P_2$ | $P_3$ | $P_4$ | $P_5$ | $P_6$ |
|---|---|---|---|---|---|
| $a$ | $a$ | $c$ | $b$ | $b$ | $c$ |
| $c$ | $c$ | $a$ | $c$ | $c$ | $b$ |
| $b$ | $b$ | $b$ | $a$ | $a$ | $a$ |

Let us calculate the value $n^{+}(x, \vec{P}_X)$ for each alternative $x \in X$: $n^{+}(a, \vec{P}_X) = 2$, $n^{+}(b, \vec{P}_X) = 2$, $n^{+}(c, \vec{P}_X) = 1$. According to the plurality rule used on the first stage the alternatives $a$ and $b$ will be chosen, i.e., $C(\vec{P}_X, X) = \{a, b\}$.

Consider now the subset $X' = X \backslash \{b\}$. A contraction of a profile $\vec{P}_X$ onto a set $X'$, i.e., $\vec{P}_{X'}$, looks as

| $P_1$ | $P_2$ | $P_3$ | $P_4$ | $P_5$ | $P_6$ |
|---|---|---|---|---|---|
| $a$ | $a$ | $c$ | $c$ | $c$ | $c$ |
| $c$ | $c$ | $a$ | $a$ | $a$ | $a$ |

According to the plurality rule the alternative $c$ will be chosen, i.e., $C(\vec{P}_{X'}, X') = \{c\}$.

Finally, let us consider the subset $X'' = X \setminus \{a\}$. A contraction of a profile $\vec{P}_X$ onto a set $X''$, i.e., $\vec{P}_{X''}$, looks as

| **P₁** | **P₂** | **P₃** | **P₄** | **P₅** | **P₅** |
|---|---|---|---|---|---|
| $c$ | $c$ | $c$ | $b$ | $b$ | $c$ |
| $b$ | $b$ | $b$ | $c$ | $c$ | $b$ |

According to the rule the alternative *c* will be chosen, i.e., $C(\vec{P}_{X''}, X'') = \{c\}$.

Then $\{a, b\} \in C(\vec{P}_X, X)$, $\{a\} \notin C(\vec{P}_{X'}, X')$ and $\{b\} \notin C(\vec{P}_{X''}, X'')$. Thus, the two-stage choice procedure does not satisfy the Monotonicity conditionon 2.

- Two-stage choice procedure 50 does not satisfy the non-compensatory condition. To prove it the following example is used.

Let $X = \{a, b, c\}$ and the profile $\vec{P}_X$ is the following

| **P₁** | **P₂** | **P₃** |
|---|---|---|
| $a$ | $a$ | $b$ |
| $c$ | $c$ | $c$ |
| $b$ | $b$ | $a$ |

Let us calculate the value $n^+(x, \vec{P}_X)$ for each alternative $x \in X$: $n^+(a, \vec{P}_X) = 2$, $n^+(b, \vec{P}_X) = 1$, $n^+(c, \vec{P}_X) = 0$. According to the rule the alternative *a* will be chosen, i.e., $C(\vec{P}_X, X) = \{a\}$.

Consider now a profile $\overrightarrow{P'_X}$, which differs from the profile $\vec{P}_X$ only by improved position of the alternative *c* in $P'_2$:

| $\mathbf{P'_1}$ | $\mathbf{P'_2}$ | $\mathbf{P'_3}$ |
|---|---|---|
| $a$ | $c$ | $b$ |
| $c$ | $a$ | $c$ |
| $b$ | $b$ | $a$ |

According to the rule the alternatives *a,b,c* will be chosen, i.e., $C(\overrightarrow{P'_X}, X) = \{a, b, c\}$.

$$C\left(\overrightarrow{P'_X}, X\right) \neq \left[\begin{array}{c} C(\vec{P}_X, X) \ or \\ \{c\} \ or \\ C(\vec{P}_X, X) \cup \{c\}. \end{array}\right.$$

Thus, the strict monotonicity condition is not satisfied.

## Two-stage choice procedures 57-84.

Two-stage choice procedures 57-84 do not satisfy the same conditions as choice procedures which are used on the second stage. To prove it we can use the same examples but with larger number of alternatives which are regarded as the worst by the maximum number of criteria [11].

Thus, it is necessary to check those normative conditions of two-stage choice procedures which were satisfied for the second stage choice procedures.

- Two-stage choice procedures 57,75,76 do not satisfy the condition **H**. To prove it the following example is used.

Let $X=\{a,b,c,d\}$ and the profile $\vec{P}_X$ is the following

| $P_1$ | $P_2$ | $P_3$ |
|---|---|---|
| $b$ | $d$ | $a$ |
| $a$ | $b$ | $d$ |
| $c$ | $a$ | $c$ |
| $d$ | $c$ | $b$ |

According to the rule the alternative *a* will be chosen, i.e., $C(\vec{P}_X,X)=\{a\}$.

Consider now the subset $X'=X\backslash\{c,d\}$. A contraction of a profile $\vec{P}_X$ onto a set $X'$, i.e., $\vec{P}_{X'}$, looks as

| $P_1$ | $P_2$ | $P_3$ |
|---|---|---|
| $b$ | $b$ | $a$ |
| $a$ | $a$ | $b$ |

According to the rule the alternative *b* will be chosen, i.e., $C(\vec{P}_{X'},X')=\{b\}$.

Then $C(\vec{P}_{X'},X')\not\supseteq C(\vec{P}_X,X)\cap X'$. Thus, the condition **H** is not satisfied.

- Two-stage choice procedures 57,68,75 do not satisfy the condition **O** (see the previous example).
- Two-stage choice procedures 68,73,75,76 do not satisfy the condition **C**. To prove it the following example is used.

Consider the previous example. According to the rule $C(\vec{P}_X,X)=\{a\}$.

Consider now the subset $X'=X\backslash\{c\}$. A contraction of a profile $\vec{P}_X$ onto a set $X'$, i.e., $\vec{P}_{X'}$, looks as

| $P_1$ | $P_2$ | $P_3$ |
|---|---|---|
| $b$ | $d$ | $a$ |
| $a$ | $b$ | $d$ |
| $d$ | $a$ | $b$ |

According to the rule the alternatives *a,b,d* will be chosen, i.e., $C(\vec{P}_{X'},X')=\{a,b,d\}$.

Finally, let us consider the subset $X''=X\backslash\{a,d\}$. A contraction of a profile $\vec{P}_X$ onto a set $X''$, i.e., $\vec{P}_{X''}$, looks as

| $P_1$ | $P_2$ | $P_3$ |
|---|---|---|
| $b$ | $b$ | $c$ |
| $c$ | $c$ | $b$ |

According to the rule the alternative *b* will be chosen, i.e., $C(\vec{P}_{X''},X'')=\{b\}$.

Then $C(\vec{P}_{X'},X')\cap C(\vec{P}_{X''},X'')=\{b\}\nsubseteq C(\vec{P}_X,X)$. Thus, the condition **C** is not satisfied.

- Two-stage choice procedures 57-60, 63,64, 68-84 satisfy the Monotonicity conditionon 1 (the proof follows from the properties of the inverse plurality rule and second stage choice procedures).
- Two-stage choice procedures 76,82 do not satisfy the Monotonicity conditionon 2. To prove it the following example is used.

Let $X=\{a,b,c\}$ and the profile $\vec{P}_X$ is the following

| $P_1$ | $P_2$ | $P_3$ | $P_4$ | $P_5$ | $P_6$ | $P_7$ | $P_8$ | $P_9$ | $P_{10}$ |
|---|---|---|---|---|---|---|---|---|---|
| $c$ | $c$ | $c$ | $a$ | $a$ | $b$ | $b$ | $c$ | $c$ | $c$ |
| $a$ | $a$ | $a$ | $b$ | $b$ | $a$ | $a$ | $b$ | $b$ | $b$ |
| $b$ | $b$ | $b$ | $c$ | $c$ | $c$ | $c$ | $a$ | $a$ | $a$ |

According to the rule the alternatives *a*,*b* will be chosen, i.e., $C(\vec{P}_X, X) = \{a, b\}$.
Consider now the subset $X' = X \backslash \{b\}$. A contraction of a profile $\vec{P}_X$ onto a set $X'$, i.e., $\vec{P}_{X'}$, looks as

| $P_1$ | $P_2$ | $P_3$ | $P_4$ | $P_5$ | $P_6$ | $P_7$ | $P_8$ | $P_9$ | $P_{10}$ |
|---|---|---|---|---|---|---|---|---|---|
| *c* | *c* | *c* | *a* | *a* | *a* | *a* | *c* | *c* | *c* |
| *a* | *a* | *a* | *c* | *c* | *c* | *c* | *a* | *a* | *a* |

According to the inverse plurality rule used on the first stage the alternative *c* will be chosen, i.e., $C(\vec{P}_{X'}, X') = \{c\}$. Finally, let us consider the subset $X'' = X \backslash \{a\}$. A contraction of a profile $\vec{P}_X$ onto a set $X''$, i.e., $\vec{P}_{X''}$, looks as

| $P_1$ | $P_2$ | $P_3$ | $P_4$ | $P_5$ | $P_6$ | $P_7$ | $P_8$ | $P_9$ | $P_{10}$ |
|---|---|---|---|---|---|---|---|---|---|
| *c* | *c* | *c* | *b* | *b* | *b* | *b* | *c* | *c* | *c* |
| *b* | *b* | *b* | *c* | *c* | *c* | *c* | *b* | *b* | *b* |

According to the rule the alternative *c* will be chosen, i.e., $C(\vec{P}_{X''}, X'') = \{c\}$.
Then $\{a, b\} \in C(\vec{P}_X, X)$, $\{a\} \notin C(\vec{P}_{X'}, X')$ and $\{b\} \notin C(\vec{P}_{X''}, X'')$. Thus, the Monotonicity conditionon 2 is not satisfied.

- Two-stage choice procedure 78 does not satisfy the non-compensatory condition. To prove it let $X = \{a, b, c\}$ and the profile $\vec{P}_X$ is the following

| $P_1$ | $P_2$ | $P_3$ | $P_4$ |
|---|---|---|---|
| *a* | *b* | *c* | *b* |
| *b* | *a* | *a* | *c* |
| *c* | *c* | *b* | *a* |

According to the non-compensatory condition the alternative *b* is better than the alternative *a* which is better than the alternative *c*. However, according to the two-stage choice procedure 78, the alternatives *a*,*b* will be chosen, i.e., $C(\vec{P}_X, X) = \{a, b\}$. Thus, the non-compensatory condition is not satisfied as $\{b\} \neq C(\vec{P}_X, X)$.

## Two-stage choice procedures 85-112.

Two-stage choice procedures 85-112 satisfy the same conditions as two-stage choice procedures 30-46 which use the plurality rule on the first stage. The proof follows from properties of q-Approval rule [11].

## Two-stage choice procedures 169-196.

Two-stage choice procedures 169-196 do not satisfy the same conditions as choice procedures which are used on the second stage. To prove it we can use the same examples [11] but with larger number of alternatives which are omitted on the first stage that were used to check the properties of the Borda rule.
Thus, it is necessary to check those normative conditions of two-stage choice procedures which were satisfied for the second stage choice procedures.

- Two-stage choice procedures 169,187,188 do not satisfy the condition **H**. To prove it the following example is provided.

Let $X = \{a, b, c, d, e\}$ and the profile $\vec{P}_X$ is the following

| $P_1$ | $P_2$ | $P_3$ | $P_4$ | $P_5$ |
|---|---|---|---|---|
| $e$ | $e$ | $a$ | $b$ | $b$ |
| $a$ | $a$ | $d$ | $c$ | $c$ |
| $b$ | $b$ | $c$ | $e$ | $a$ |
| $c$ | $c$ | $b$ | $a$ | $d$ |
| $d$ | $d$ | $e$ | $d$ | $e$ |

Let us calculate the Borda count for each alternative: $r(a, \vec{P}_X) = r(b, \vec{P}_X) = 13, r(c, \vec{P}_X) = 10$, $r(d, \vec{P}_X) = 4, r(e, \vec{P}_X) = 10$. According to the rule the alternatives *c*,*d*,*e* will be omitted on the first stage and the alternative *a* will be chosen on the second stage, i.e., $C(\vec{P}_X, X) = \{a\}$.

Consider now the subset $X' = X \backslash \{d\}$. A contraction of a profile $\vec{P}_X$ onto a set $X'$, i.e., $\vec{P}_{X'}$, looks as

| $P_1$ | $P_2$ | $P_3$ | $P_4$ | $P_5$ |
|---|---|---|---|---|
| $e$ | $e$ | $a$ | $b$ | $b$ |
| $a$ | $a$ | $c$ | $c$ | $c$ |
| $b$ | $b$ | $b$ | $e$ | $a$ |
| $c$ | $c$ | $e$ | $a$ | $e$ |

Let us calculate the Borda count for each alternative: $r(a, \vec{P}_{X'}) = 8, r(b, \vec{P}_{X'}) = 9, r(c, \vec{P}_{X'}) = 6$, $r(e, \vec{P}_{X'}) = 7$. According to the rule the alternative *b* will be chosen, i.e., $C(\vec{P}_{X'}, X') = \{b\}$.

Then $C(\vec{P}_{X'}, X') \not\supseteq C(\vec{P}_X, X) \cap X'$. Thus, the condition **H** is not satisfied.

- Two-stage choice procedures 169, 180, 187 do not satisfy the condition **O** (see the previous example).
- Two-stage choice procedures 180,185,187,188 do not satisfy the condition **C** (see paragraph 7.2 of the previous section).

Let $X = \{a, b, c, d\}$ and the profile $\vec{P}_X$ is the following

| $P_1$ | $P_2$ | $P_3$ |
|---|---|---|
| $a$ | $b$ | $c$ |
| $c$ | $a$ | $b$ |
| $d$ | $d$ | $a$ |
| $b$ | $c$ | $d$ |

Let us calculate the Borda count for each alternative: $r(a, \vec{P}_X) = 6, r(b, \vec{P}_X) = r(c, \vec{P}_X) = 5$, $r(d, \vec{P}_X) = 2$. According to the rule the alternative *a* will be chosen, i.e., $C(\vec{P}_X, X) = \{a\}$.

Consider now the subset $X' = X \backslash \{d\}$. A contraction of a profile $\vec{P}_X$ onto a set $X'$, i.e., $\vec{P}_{X'}$, looks as

| $P_1$ | $P_2$ | $P_3$ |
|---|---|---|
| $a$ | $b$ | $c$ |
| $c$ | $a$ | $b$ |
| $b$ | $c$ | $a$ |

Let us calculate the Borda count for each alternative: $r(a, \vec{P}_{X'}) = r(b, \vec{P}_{X'}) = r(c, \vec{P}_{X'}) = 3$. According to the rule the alternatives *a*,*b*,*c* will be chosen, i.e., $C(\vec{P}_{X'}, X') = \{a, b, c\}$.

Finally, let us consider the subset $X'' = X \backslash \{a\}$. A contraction of a profile $\vec{P}_X$ onto a set $X''$, i.e., $\vec{P}_{X''}$, looks as

| $\mathbf{P_1}$ | $\mathbf{P_2}$ | $\mathbf{P_3}$ |
|---|---|---|
| *c* | *b* | *c* |
| *d* | *d* | *b* |
| *b* | *c* | *d* |

Let us calculate the Borda count for each alternative: $r(b, \vec{P}_{X''}) = 3, r(c, \vec{P}_{X''}) = 4, r(d, \vec{P}_{X''}) = 2$.
According to the rule the alternative *c* will be chosen, i.e., $C(\vec{P}_{X''}, X'') = \{c\}$.
Then $C(\vec{P}_{X'}, X') \cap C(\vec{P}_{X''}, X'') = \{c\} \not\subseteq C(\vec{P}_X, X)$. Thus, the condition **C** is not satisfied.

- Two-stage choice procedures 169-196 satisfy the Monotonicity conditionon 1 (the proof follows from properties of the Borda rule).
- Two-stage choice procedures 188, 194 do not satisfy the Monotonicity conditionon 2. To prove it the following example is provided.

Let $X = \{a, b, c, d, e, f\}$ and the profile $\vec{P}_X$ is the following

| $\mathbf{P_1}$ | $\mathbf{P_2}$ | $\mathbf{P_3}$ | $\mathbf{P_4}$ | $\mathbf{P_5}$ | $\mathbf{P_6}$ | $\mathbf{P_7}$ | $\mathbf{P_8}$ |
|---|---|---|---|---|---|---|---|
| *a* | *a* | *a* | *f* | *f* | *f* | *f* | *c* |
| *b* | *b* | *b* | *b* | *b* | *b* | *c* | *f* |
| *c* | *c* | *c* | *a* | *a* | *a* | *e* | *e* |
| *d* | *d* | *d* | *c* | *c* | *c* | *d* | *d* |
| *e* | *e* | *e* | *d* | *d* | *d* | *b* | *a* |
| *f* | *f* | *f* | *e* | *e* | *e* | *a* | *b* |

Let us calculate the Borda count for each alternative: $r(a, \vec{P}_X) = r(b, \vec{P}_X) = 25$, $r(c, \vec{P}_X) = 24$, $r(d, \vec{P}_X) = 13, r(e, \vec{P}_X) = 9, r(f, \vec{P}_X) = 24$. According to the rule the alternatives *a*,*b* will be chosen, i.e., $C(\vec{P}_X, X) = \{a, b\}$.
Consider now the subset $X' = X \backslash \{b\}$. A contraction of a profile $\vec{P}_X$ onto a set $X'$, i.e., $\vec{P}_{X'}$, looks as

| $\mathbf{P_1}$ | $\mathbf{P_2}$ | $\mathbf{P_3}$ | $\mathbf{P_4}$ | $\mathbf{P_5}$ | $\mathbf{P_6}$ | $\mathbf{P_7}$ | $\mathbf{P_8}$ |
|---|---|---|---|---|---|---|---|
| *a* | *a* | *a* | *f* | *f* | *f* | *f* | *c* |
| *c* | *c* | *c* | *a* | *a* | *a* | *c* | *f* |
| *d* | *d* | *d* | *c* | *c* | *c* | *e* | *e* |
| *e* | *e* | *e* | *d* | *d* | *d* | *d* | *d* |
| *f* | *f* | *f* | *e* | *e* | *e* | *a* | *a* |

Let us calculate the Borda count for each alternative: $r(a, \vec{P}_{X'}) = 21, r(c, \vec{P}_{X'}) = 22$, $r(d, \vec{P}_{X'}) = 13$, $r(e, \vec{P}_{X'}) = 7, r(f, \vec{P}_{X'}) = 19$. According to the rule the alternative *c* will be chosen, i.e., $C(\vec{P}_{X'}, X') = \{c\}$.
Finally, let us consider the subset $X'' = X \backslash \{a\}$. A contraction of a profile $\vec{P}_X$ onto a set $X''$, i.e., $\vec{P}_{X''}$, looks as

| $\mathbf{P_1}$ | $\mathbf{P_2}$ | $\mathbf{P_3}$ | $\mathbf{P_4}$ | $\mathbf{P_5}$ | $\mathbf{P_6}$ | $\mathbf{P_7}$ | $\mathbf{P_8}$ |
|---|---|---|---|---|---|---|---|
| *b* | *b* | *b* | *f* | *f* | *f* | *f* | *c* |
| *c* | *c* | *c* | *b* | *b* | *b* | *c* | *f* |
| *d* | *d* | *d* | *c* | *c* | *c* | *e* | *e* |
| *e* | *e* | *e* | *d* | *d* | *d* | *d* | *d* |
| *f* | *f* | *f* | *e* | *e* | *e* | *b* | *b* |

Let us calculate the Borda count for each alternative: $r(b,\vec{P}_{X''}) = 21, r(c,\vec{P}_{X''}) = 22, r(d,\vec{P}_{X''}) = 13$, $r(e,\vec{P}_{X''}) = 7, r(f,\vec{P}_{X''}) = 19$. According to the rule the alternative $c$ will be chosen, i.e., $C(\vec{P}_{X''},X'') = \{c\}$.

Then $\{a,b\} \in C(\vec{P}_X,X)$, $\{a\} \notin C(\vec{P}_{X'},X')$ and $\{b\} \notin C(\vec{P}_{X''},X'')$. Thus, the Monotonicity conditionon 2 is not satisfied.

- Two-stage choice procedure 190 does not satisfy the non-compensatory condition. To prove it let $X = \{a,b,c\}$ and the profile $\vec{P}_X$ is the following

| **P1** | **P2** | **P3** |
|---|---|---|
| *a* | *a* | *b* |
| *b* | *b* | *c* |
| *c* | *c* | *a* |

According to the two-stage choice procedure the alternative *a* will be chosen, i.e., $C(\vec{P}_X,X) = \{a\}$. However, according to the non-compensatory condition the alternative *b* is better than the alternative *a* which is better than the alternative *c*.

Thus, the non-compensatory condition is not satisfied as $\{b\} \neq C(\vec{P}_X,X)$.

## Two-stage choice procedures 197-224.

To check the properties of two-stage choice procedures 197-224 similar examples, which were given to check the properties of existing choice procedures [11], can be used.

- Two-stage choice procedures 197-224 do not satisfy the conditions **H** and **O.** To prove it the following example is provided.

Let $X = \{a,b,c\}$ and the profile $\vec{P}_X$ is the following

| **P1** | **P2** | **P3** | **P4** | **P5** |
|---|---|---|---|---|
| *a* | *a* | *d* | *b* | *b* |
| *b* | *b* | *a* | *d* | *d* |
| *d* | *d* | *b* | *a* | *a* |

In this example none of alternatives is a Condorcet winner ($a\mu b, d\mu a, b\mu d$). Thus, the Borda rule is used for this case. Let us calculate the Borda count for each alternative: $r(a,\vec{P}_X) = 5, r(b,\vec{P}_X) = 6, r(c,\vec{P}_X) = 4$. According to the rule the alternative *b* will be chosen, i.e., $C(\vec{P}_X,X) = \{b\}$.

Consider now the subset $X' = X\backslash\{d\}$. A contraction of a profile $\vec{P}_X$ onto a set $X'$, i.e., $\vec{P}_{X'}$, looks as

| **P1** | **P2** | **P3** | **P4** | **P5** |
|---|---|---|---|---|
| *a* | *a* | *a* | *b* | *b* |
| *b* | *b* | *b* | *a* | *a* |

According to the rule the alternative *a* will be chosen, i.e., $C(\vec{P}_{X'},X') = \{a\}$.

$C(\vec{P}_{X'},X') \not\supseteq C(\vec{P}_X,X) \cap X'$. Thus, the condition **H** is not satisfied.

$C(\vec{P}_{X'},X') \neq C(\vec{P}_X,X)$. Thus, the condition **O** is not satisfied.

- Two-stage choice procedures 197-224 do not satisfy the condition **C**. To prove it the following example is provided.

Let $X = \{a, b, c, d\}$ and the profile $\vec{P}_X$ is the following

| **P₁** | **P₂** | **P₃** | **P₄** | **P₅** |
|---|---|---|---|---|
| *a* | *a* | *c* | *d* | *d* |
| *d* | *d* | *b* | *b* | *b* |
| *c* | *c* | *a* | *a* | *a* |
| *b* | *b* | *d* | *c* | *c* |

In this example none of alternatives is a Condorcet winner ($a\mu c, a\mu d, b\mu a, c\mu b, d\mu b, d\mu c$). Thus, the Borda rule is used for this case. Let us calculate the Borda count for each alternative: $r(a, \vec{P}_X) = 9, r(b, \vec{P}_X) = 6, r(c, \vec{P}_X) = 5$, $r(d, \vec{P}_X) = 10$. According to the rule the alternative *d* will be chosen, i.e., $C(\vec{P}_X, X) = \{d\}$.

Consider now the subset $X' = X \backslash \{d\}$. A contraction of a profile $\vec{P}_X$ onto a set $X'$, i.e., $\vec{P}_{X'}$, looks as

| **P₁** | **P₂** | **P₃** | **P₄** | **P₅** |
|---|---|---|---|---|
| *a* | *a* | *c* | *b* | *b* |
| *c* | *c* | *b* | *a* | *a* |
| *b* | *b* | *a* | *c* | *c* |

According to the rule the alternative *a* will be chosen, i.e., $C(\vec{P}_{X'}, X') = \{a\}$.

Finally, let us consider the subset $X'' = X \backslash \{b, c\}$. A contraction of a profile $\vec{P}_X$ onto a set $X''$, i.e., $\vec{P}_{X''}$, looks as

| **P₁** | **P₂** | **P₃** | **P₄** | **P₅** |
|---|---|---|---|---|
| *a* | *a* | *a* | *d* | *d* |
| *d* | *d* | *d* | *a* | *a* |

According to the rule the alternative *a* will be chosen, i.e., $C(\vec{P}_{X''}, X'') = \{a\}$.

Then $C(\vec{P}_{X'}, X') \cap C(\vec{P}_{X''}, X'') = \{a\} \nsubseteq C(\vec{P}_X, X)$. Thus, the condition **C** is not satisfied.

- Two-stage choice procedures 197-224 do not satisfy the condition **ACA** since such choice procedures do not satisfy the condition **H**.
- Two-stage choice procedures 197-224 satisfy the Monotonicity conditionon 1 (the proof follows from the properties of the Black procedure).
- Two-stage choice procedures 198-200, 203-214, 216-224 do not satisfy the Monotonicity conditionon 2 (the proof follows from the properties of the Black procedure).
- Since two-stage choice procedures 197, 201-202, 215 chooses no more than one best alternative, the Monotonicity conditionon 2 is not applicable to these procedures as it considers the choice of more than two alternatives. In other words, such two-stage choice procedures obey the Monotonicity conditionon 2 trivially.
- Two-stage choice procedures 197-224 do not satisfy the strict monotonicity condition (the proof follows from the properties of the Black procedure).
- Two-stage choice procedures 197-224 do not satisfy the non-compensatory condition (the proof follows from the properties of the Black procedure).

## Two-stage choice procedures 225-252.

Two-stage choice procedures 225-252 do not satisfy the same properties as the Inverse Borda rule. To prove it we can use the same examples that were used to check the properties of the Inverse Borda rule [11].

### Two-stage choice procedures 253-280.

Two-stage choice procedures 253-280 do not satisfy the same properties as the Nanson rule. To prove it we can use the same examples that were used to check the properties of the Nanson rule [11].

### Two-stage choice procedures 281-308.

Two-stage choice procedures 281-308 do not satisfy the same properties as the Coombs procedure. To prove it we can use the same examples that were used to check the properties of the Coombs procedure [11].

### Two-stage choice procedures 310-319, 330, 334-336

***310.5 Monotonicity condition 1***

Let $X = \{a, b, c, d\}$ and the profile $\vec{P}_X$ is the following

| $P_1$ | $P_2$ | $P_3$ | $P_4$ | $P_5$ | $P_6$ | $P_7$ | $P_8$ | $P_9$ | $P_{10}$ |
|---|---|---|---|---|---|---|---|---|---|
| *b* | *d* | *a* | *a* | *b* | *d* | *b* | *d* | *a* | *c* |
| *d* | *c* | *b* | *c* | *a* | *c* | *c* | *c* | *b* | *a* |
| *a* | *a* | *d* | *b* | *c* | *b* | *d* | *a* | *d* | *d* |
| *c* | *b* | *c* | *d* | *d* | *a* | *a* | *b* | *c* | *b* |

According to the rule the alternatives *a*,*b*,*d* will be chosen, i.e., $C(\vec{P}_X, X) = \{a, b, d\}$.

Consider now a profile $\overrightarrow{P'_X}$, which differs from the profile $\vec{P}_X$ only by improved position of the alternative *a* in $P'_7$:

| $P'_1$ | $P'_2$ | $P'_3$ | $P'_4$ | $P'_5$ | $P'_6$ | $P'_7$ | $P'_8$ | $P'_9$ | $P'_{10}$ |
|---|---|---|---|---|---|---|---|---|---|
| *b* | *d* | *a* | *a* | *b* | *d* | *b* | *d* | *a* | *c* |
| *d* | *c* | *b* | *c* | *a* | *c* | *c* | *c* | *b* | *a* |
| *a* | *a* | *d* | *b* | *c* | *b* | *a* | *a* | *d* | *d* |
| *c* | *b* | *c* | *d* | *d* | *a* | *d* | *b* | *c* | *b* |

According to the two-stage choice procedure the alternative *d* is omitted on the first stage. Consider now the subset $X' = X\backslash\{d\}$. A contraction of a profile $\overrightarrow{P'_X}$ onto a set $X'$, i.e., $\overrightarrow{P'_{X'}}$, looks as

| $P'_1$ | $P'_2$ | $P'_3$ | $P'_4$ | $P'_5$ | $P'_6$ | $P'_7$ | $P'_8$ | $P'_9$ | $P'_{10}$ |
|---|---|---|---|---|---|---|---|---|---|
| *b* | *c* | *a* | *a* | *b* | *c* | *b* | *c* | *a* | *c* |
| *a* | *a* | *b* | *c* | *a* | *b* | *c* | *a* | *b* | *a* |
| *c* | *b* | *c* | *b* | *c* | *a* | *a* | *b* | *c* | *b* |

According to the rule the alternative *c* will be chosen, i.e., $C\left(\overrightarrow{P'_X}, X\right) = \{c\}$.

Then $\{a\} \in C(\vec{P}_X, X)$, $\{a\} \notin C\left(\overrightarrow{P'_X}, X\right)$. Thus, the Monotonicity conditionon 1 is not satisfied.

***311.5 Monotonicity condition 1***

Let $X = \{a, b, c, d\}$ and the profile $\vec{P}_X$ is the following

| $P_1$ | $P_2$ | $P_3$ | $P_4$ | $P_5$ | $P_6$ | $P_7$ | $P_8$ | $P_9$ | $P_{10}$ |
|---|---|---|---|---|---|---|---|---|---|
| *b* | *d* | *c* | *a* | *a* | *d* | *b* | *d* | *a* | *d* |
| *d* | *a* | *b* | *c* | *b* | *c* | *c* | *c* | *b* | *a* |
| *c* | *c* | *a* | *b* | *c* | *b* | *a* | *a* | *d* | *b* |
| *a* | *b* | *d* | *d* | *d* | *a* | *d* | *b* | *c* | *c* |

According to the rule the alternatives *a*,*b*,*c* will be chosen, i.e., $C(\vec{P}_X, X) = \{a, b, c\}$.

Consider now a profile $\overrightarrow{P'_X}$, which differs from the profile $\vec{P}_X$ only by improved position of the alternative *a* in $P'_2$:

| $P'_1$ | $P'_2$ | $P'_3$ | $P'_4$ | $P'_5$ | $P'_6$ | $P'_7$ | $P'_8$ | $P'_9$ | $P'_{10}$ |
|---|---|---|---|---|---|---|---|---|---|
| $b$ | $a$ | $c$ | $a$ | $a$ | $d$ | $b$ | $d$ | $a$ | $d$ |
| $d$ | $d$ | $b$ | $c$ | $b$ | $c$ | $c$ | $c$ | $b$ | $a$ |
| $c$ | $c$ | $a$ | $b$ | $c$ | $b$ | $a$ | $a$ | $d$ | $b$ |
| $a$ | $b$ | $d$ | $d$ | $d$ | $a$ | $d$ | $b$ | $c$ | $c$ |

According to the two-stage choice procedure the alternative *d* is omitted on the first stage. Consider now the subset $X' = X\backslash\{d\}$. A contraction of a profile $\overrightarrow{P'_X}$ onto a set $X'$, i.e., $\overrightarrow{P'_{X'}}$, looks as

| $P'_1$ | $P'_2$ | $P'_3$ | $P'_4$ | $P'_5$ | $P'_6$ | $P'_7$ | $P'_8$ | $P'_9$ | $P'_{10}$ |
|---|---|---|---|---|---|---|---|---|---|
| $b$ | $a$ | $c$ | $a$ | $a$ | $c$ | $b$ | $c$ | $a$ | $a$ |
| $c$ | $c$ | $b$ | $c$ | $b$ | $b$ | $c$ | $a$ | $b$ | $b$ |
| $a$ | $b$ | $a$ | $b$ | $c$ | $a$ | $a$ | $b$ | $c$ | $c$ |

According to the rule the alternatives *b*,*c* will be chosen, i.e., $C\left(\overrightarrow{P'_X}, X\right) = \{b, c\}$.

Then $\{a\} \in C(\vec{P}_X, X)$, $\{a\} \notin C\left(\overrightarrow{P'_X}, X\right)$. Thus, the Monotonicity conditionon 1 is not satisfied.

### *315.5 Monotonicity condition 1*

Let $X = \{a, b, c, d\}$ and the profile $\vec{P}_X$ is the following

| $P_1$ | $P_2$ | $P_3$ | $P_4$ | $P_5$ | $P_6$ | $P_7$ | $P_8$ | $P_9$ | $P_{10}$ |
|---|---|---|---|---|---|---|---|---|---|
| $b$ | $d$ | $c$ | $a$ | $a$ | $d$ | $b$ | $d$ | $a$ | $d$ |
| $d$ | $c$ | $b$ | $c$ | $b$ | $c$ | $c$ | $c$ | $b$ | $a$ |
| $c$ | $a$ | $a$ | $b$ | $c$ | $b$ | $a$ | $a$ | $d$ | $b$ |
| $a$ | $b$ | $d$ | $d$ | $d$ | $a$ | $d$ | $b$ | $c$ | $c$ |

Let us calculate the Borda count for each alternative: $r(a, \vec{P}_X) = r(b, \vec{P}_X) = r(c, \vec{P}_X) = r(d, \vec{P}_X) = 15$. According to the rule the alternatives *a*,*b*,*c*,*d* will be chosen, i.e., $C(\vec{P}_X, X) = \{a, b, c, d\}$.
Consider now a profile $\overrightarrow{P'_X}$, which differs from the profile $\vec{P}_X$ only by improved position of the alternative *a* in $P'_{10}$:

| $P'_1$ | $P'_2$ | $P'_3$ | $P'_4$ | $P'_5$ | $P'_6$ | $P'_7$ | $P'_8$ | $P'_9$ | $P'_{10}$ |
|---|---|---|---|---|---|---|---|---|---|
| $b$ | $d$ | $c$ | $a$ | $a$ | $d$ | $b$ | $d$ | $a$ | $a$ |
| $d$ | $c$ | $b$ | $c$ | $b$ | $c$ | $c$ | $c$ | $b$ | $d$ |
| $c$ | $a$ | $a$ | $b$ | $c$ | $b$ | $a$ | $a$ | $d$ | $b$ |
| $a$ | $b$ | $d$ | $d$ | $d$ | $a$ | $d$ | $b$ | $c$ | $c$ |

According to the two-stage choice procedure the alternative *d* is omitted on the first stage. Consider now the subset $X' = X\backslash\{d\}$. A contraction of a profile $\overrightarrow{P'_X}$ onto a set $X'$, i.e., $\overrightarrow{P'_{X'}}$, looks as

| $P'_1$ | $P'_2$ | $P'_3$ | $P'_4$ | $P'_5$ | $P'_6$ | $P'_7$ | $P'_8$ | $P'_9$ | $P'_{10}$ |
|---|---|---|---|---|---|---|---|---|---|
| $b$ | $c$ | $c$ | $a$ | $a$ | $c$ | $b$ | $c$ | $a$ | $a$ |
| $c$ | $a$ | $b$ | $c$ | $b$ | $b$ | $c$ | $a$ | $b$ | $b$ |
| $a$ | $b$ | $a$ | $b$ | $c$ | $a$ | $a$ | $b$ | $c$ | $c$ |

According to the rule the alternative *c* will be chosen, i.e., $C\left(\overrightarrow{P'_X}, X\right) = \{c\}$.

Then $\{a\} \in C(\vec{P}_X, X)$, $\{a\} \notin C\left(\overrightarrow{P'_X}, X\right)$. Thus, the Monotonicity conditionon 1 is not satisfied.

To check the remaining normative conditions of given two-stage choice procedures similar examples which were given to check the properties of existing choice procedures can be used.

- Two-stage choice procedures 310-319, 330, 334-336 do not satisfy the condition **H** (see the properties of minimal dominant set).

- Two-stage choice procedures 310-319, 330, 334-336 do not satisfy the condition **C** (see the properties of the procedures used on the second stage).
- Two-stage choice procedures 310-319, 330, 334-336 do not satisfy the condition **O** (see the properties of the procedures used on the second stage and paragraphs 310.5, 315.5 of Appendix 1).
- Two-stage choice procedures 310-319, 330, 334-336 do not satisfy the condition **ACA** since such choice procedures do not satisfy the condition **H**.
- Two-stage choice procedures 310-319, 330 do not satisfy the Monotonicity conditionon 1 (see the properties of the Coombs procedure for choice procedure 319 and paragraphs 310.5, 311.5, 315.5 for other choice procedures).
- Two-stage choice procedures 334-336 satisfy the Monotonicity conditionon 1 (the proof follows from the properties of minimal dominant set and second stage choice procedures).
- Two-stage choice procedures 310-313, 316-319, 330, 334-336 do not satisfy the Monotonicity conditionon 2 (see the properties of minimal dominant set).
- Since two-stage choice procedures 314-315 chooses no more than one best alternative, the Monotonicity conditionon 2 is not applicable to these procedures as it considers the choice of more than two alternatives. In other words, such two-stage choice procedures obey the Monotonicity conditionon 2 trivially.
- Two-stage choice procedures 310-319, 330, 334-336 do not satisfy the strict monotonicity condition (see the properties of minimal dominant set).
- Two-stage choice procedures 310-319, 330, 334-336 do not satisfy the non-compensatory condition (see the properties of minimal dominant set).

### Two-stage choice procedures 338-347, 350-364

Since minimal undominated set is equal to minimal dominant set when the number of criteria is odd, two-stage choice procedures 338-347, 358, 362-364 do not satisfy the same conditions as two-stage choice procedures 310-319, 330, 334-336, which use minimal dominant set on the first stage.
Consider now the properties of other two-stage choice procedures.

- Two-stage choice procedures 351-354, 357, 359-361 do not satisfy the condition **H, C** and **O**. To prove it let $X = \{a, b, c, d\}$ and the profile $\vec{P}_X$ looks as

| $P_1$ | $P_2$ | $P_3$ | $P_4$ | $P_5$ | $P_6$ |
|---|---|---|---|---|---|
| *b* | *a* | *a* | *b* | *a* | *d* |
| *a* | *c* | *d* | *a* | *c* | *b* |
| *c* | *d* | *b* | *c* | *d* | *c* |
| *d* | *b* | *c* | *d* | *b* | *a* |

For this case a matrix of majority relation μ is the following

| | **a** | **b** | **c** | **d** |
|---|---|---|---|---|
| **a** | - | 0 | 1 | 1 |
| **b** | 0 | - | 1 | 0 |
| **c** | 0 | 0 | - | 1 |
| **d** | 0 | 1 | 0 | - |

According to the rule the alternative $a$ is included in minimal undominated set $Q$. Thus, $C(\vec{P}_X, X) = Q = \{a\}$.

Consider now the subset $X' = X\backslash\{a\}$. A contraction of a profile $\vec{P}_X$ onto a set $X'$, i.e., $\vec{P}_{X'}$, looks as

| **P$_1$** | **P$_2$** | **P$_3$** | **P$_4$** | **P$_5$** | **P$_6$** |
|---|---|---|---|---|---|
| $b$ | $c$ | $d$ | $b$ | $c$ | $d$ |
| $c$ | $d$ | $b$ | $c$ | $d$ | $b$ |
| $d$ | $b$ | $c$ | $d$ | $b$ | $c$ |

According to the rule the alternatives *b,c,d* make the minimal undominated set. Thus, $C(\vec{P}_{X'}, X') = \{b, c, d\}$.

Finally, consider the subset $X'' = X\backslash\{c, d\}$. A contraction of a profile $\vec{P}_X$ onto a set $X''$, i.e., $\vec{P}_{X''}$, looks as

| **P$_1$** | **P$_2$** | **P$_3$** | **P$_4$** | **P$_5$** | **P$_6$** |
|---|---|---|---|---|---|
| $b$ | $a$ | $a$ | $b$ | $a$ | $b$ |
| $a$ | $b$ | $b$ | $a$ | $b$ | $a$ |

According to the rule the alternatives *a* and *b* will be chosen, i.e., $C(\vec{P}_{X''}, X'') = \{a, b\}$.

$C(\vec{P}_{X'}, X') \not\supseteq C(\vec{P}_X, X) \cap X'$. Thus, the condition **H** is not satisfied.

$C(\vec{P}_{X'}, X') \cap C(\vec{P}_{X''}, X'') = \{b\} \not\subseteq C(\vec{P}_X, X)$. Thus, the condition **C** is not satisfied.

$C(\vec{P}_{X'}, X') \neq C(\vec{P}_X, X)$. Thus, the condition **O** is not satisfied.

- Two-stage choice procedures 351-354, 357, 359-361 do not satisfy the condition **ACA** since such two-stage choice procedures do not satisfy the condition **H**.
- Two-stage choice procedures 351-354, 357, 359-361 satisfy the Monotonicity conditionon 1 (the proof follows from properties of minimal undominated set and second stage choice procedures).
- Two-stage choice procedures 351-354, 357, 359-361 do not satisfy the Monotonicity conditionon 2 (see the properties of minimal undominated set).
- Two-stage choice procedures 351-354, 357, 359-361 do not satisfy the strict monotonicity condition (see the properties of minimal undominated set).
- Two-stage choice procedures 351-354, 357, 359-361 do not satisfy the non-compensatory condition (see the properties of minimal undominated set).

## Two-stage choice procedures 365-392

Two-stage choice procedures do not satisfy the same conditions as choice procedures used on the second stage. To prove it we can use the same examples but with larger number of alternatives eliminated on the first stage of the two-stage choice procedure [11].

Thus, it is necessary to check those normative conditions of two-stage choice procedures which were satisfied for the second stage choice procedures.

- Two-stage choice procedures 365,383,384 do not satisfy the condition **H** (see paragraph 365.1 of Appendix 1).
- Two-stage choice procedures 365, 376 do not satisfy the condition **O** (see paragraph 365.1 of Appendix 1).
- Two-stage choice procedures 376, 381, 383, 384 do not satisfy the condition **C**. To prove it the following example is used.

Let $X = \{a, b, c, d, e\}$ and a matrix of majority relation μ is the following

| | a | b | c | d | e |
|---|---|---|---|---|---|
| **a** | - | 1 | 0 | 0 | 0 |
| **b** | 0 | - | 1 | 0 | 0 |
| **c** | 0 | 0 | - | 1 | 1 |
| **d** | 1 | 0 | 0 | - | 0 |
| **e** | 0 | 0 | 0 | 1 | - |

According to the rule the alternatives $a$ and $c$ are included in minimal weakly stable set $Q$. Thus, $C(\vec{P}_X, X) = Q = \{a, c\}$.

Consider now the subset $X' = X \backslash \{e\}$. Then a matrix of majority relation μ is the following

| | a | b | c | d |
|---|---|---|---|---|
| **a** | - | 1 | 0 | 0 |
| **b** | 0 | - | 1 | 0 |
| **c** | 0 | 0 | - | 1 |
| **d** | 1 | 0 | 0 | - |

According to the rule the alternatives *a,b,c,d* are included in minimal weakly stable set, i.e., $C(\vec{P}_{X'}, X') = \{a, b, c, d\}$.

Finally, consider the subset $X'' = X \backslash \{a\}$. Then a matrix of majority relation μ is the following

| | b | c | d | e |
|---|---|---|---|---|
| **b** | - | 1 | 0 | 0 |
| **c** | 0 | - | 1 | 1 |
| **d** | 0 | 0 | - | 0 |
| **e** | 0 | 0 | 1 | - |

According to the rule the alternative $b$ is included in minimal weakly stable set, i.e., $C(\vec{P}_{X''}, X'') = \{b\}$.

Then $C(\vec{P}_{X'}, X') \cap C(\vec{P}_{X''}, X'') = \{b\} \nsubseteq C(\vec{P}_X, X)$. Thus, the condition **C** is not satisfied.

- Two-stage choice procedures 383, 384 do not satisfy the condition **C** (see paragraph 383.2 of Appendix 1).
- Two-stage choice procedures 365-375, 383-384, 386, 390-392 do not satisfy the Monotonicity conditionon 1 (the proof follows from the paragraph 365.5 of Appendix 1).
- Two-stage choice procedures 376-382, 385, 387-389 do not satisfy the Monotonicity conditionon 1 (the proof follows from the paragraph 376.5 of Appendix 1).
- Two-stage choice procedures 365-392 do not satisfy the Monotonicity conditionon 2 (see the properties of minimal weakly stable set).
- Two-stage choice procedures 365-392 do not satisfy the non-compensatory condition (see the properties of minimal weakly stable set).

***365.1 Heredity condition (H)***

Let $X = \{a, b, c, d, e, f\}$ and a matrix of majority relation μ is the following

| | **a** | **b** | **c** | **d** | **e** | **f** |
|---|---|---|---|---|---|---|
| **a** | - | 1 | 0 | 0 | 0 | 1 |
| **b** | 0 | - | 1 | 0 | 0 | 0 |
| **c** | 1 | 0 | - | 0 | 0 | 0 |
| **d** | 0 | 1 | 0 | - | 0 | 0 |
| **e** | 0 | 0 | 0 | 1 | - | 0 |
| **f** | 0 | 0 | 0 | 0 | 1 | - |

Minimal weakly stable set is chosen on the first stage of choice procedure. According to the rule the alternatives *a*,*c* are included in minimal weakly stable set *Q*, i.e., $Q = \{a, c\}$.

Simple majority rule is used on the second stage of the choice procedure. Consider now the remaining set $X' = X \backslash \{b, d, e, f\}$. A contraction of a profile $\vec{P}_X$ onto a set $X'$, i.e., $\vec{P}_{X'}$, looks as

| $P_1$ | $P_2$ | $P_3$ |
|---|---|---|
| *c* | *c* | *a* |
| *a* | *a* | *c* |

According to the simple majority rule, the alternative *c* will be chosen, i.e., $C(\vec{P}_X, X) = \{c\}$.

Consider now the subset $X' = X \backslash \{d, e, f\}$. Then, according to the two-stage choice procedure, the alternatives *a*,*b*,*c* are included in minimal weakly stable set $Q'$.

Simple majority rule is used on the second stage of the choice procedure. Consider now the remaining set $Q' = X \backslash \{b, d, e, f\}$. A contraction of a profile $\vec{P}_X$ onto a set $Q'$, i.e., $\vec{P}_{Q'}$, looks as

| $P_1$ | $P_2$ | $P_3$ |
|---|---|---|
| *c* | *b* | *a* |
| *a* | *c* | *b* |
| *b* | *a* | *c* |

According to the rule the choice is empty, i.e., $C(\vec{P}_{X'}, X') = \emptyset$.

Then $C(\vec{P}_X, X) \cap X' = \{c\} \nsubseteq C(\vec{P}_{X'}, X')$. Thus, the condition **H** is not satisfied.

***365.5 Monotonicity condition 1.***

Consider the example from paragraph 365.1. According to the rule the alternative *c* will be chosen, i.e., $C(\vec{P}_X, X) = \{\text{c}\}$.

Suppose that the relative comparison of the alternative *c* and the alternative *d* is changed (*cμd*) while the relative comparison of any pair of other alternatives remains unchanged. Then a matrix of majority relation looks as

| | **a** | **b** | **c** | **d** | **e** | **f** |
|---|---|---|---|---|---|---|
| **a** | - | 1 | 0 | 0 | 0 | 1 |
| **b** | 0 | - | 1 | 0 | 0 | 0 |
| **c** | 1 | 0 | - | 1 | 0 | 0 |
| **d** | 0 | 1 | 0 | - | 0 | 0 |
| **e** | 0 | 0 | 0 | 1 | - | 0 |
| **f** | 0 | 0 | 0 | 0 | 1 | - |

Minimal weakly stable set is used on the first stage of choice procedure. For this case the alternatives *a*,*b*,*c* form a minimal weakly stable set. Simple majority rule is applied on the second stage of the two-stage choice procedure. Since *cμa*, *aμb*, *bμc*, $C\left(\overrightarrow{P'_X}, X\right) = \emptyset$.

Then $\{c\} \in C(\vec{P}_X, X)$, $\{c\} \notin C\left(\overrightarrow{P'_X}, X\right)$. Thus, the Monotonicity conditionon 1 is not satisfied.

***376.5 Monotonicity condition 1.***

Let $X = \{a, b, c, d, e, f\}$ and a matrix of majority relation μ is the following

| | **a** | **b** | **c** | **d** | **e** | **f** | **g** | **h** | **i** | **j** | **k** | **l** | **m** | **n** | **o** |
|---|---|---|---|---|---|---|---|---|---|---|---|---|---|---|---|
| **a** | - | 0 | 1 | 0 | 0 | 0 | 1 | 1 | 0 | 0 | 0 | 0 | 0 | 0 | 0 |
| **b** | 1 | - | 0 | 0 | 0 | 0 | 0 | 0 | 0 | 0 | 1 | 1 | 0 | 0 | 0 |
| **c** | 0 | 1 | - | 1 | 0 | 0 | 0 | 0 | 0 | 0 | 0 | 0 | 0 | 0 | 1 |
| **d** | 1 | 0 | 0 | - | 0 | 0 | 0 | 0 | 0 | 0 | 0 | 0 | 0 | 0 | 0 |
| **e** | 0 | 0 | 0 | 1 | - | 0 | 1 | 0 | 0 | 0 | 0 | 0 | 0 | 0 | 0 |
| **f** | 0 | 0 | 0 | 0 | 1 | - | 0 | 0 | 0 | 0 | 0 | 0 | 0 | 0 | 0 |
| **g** | 0 | 0 | 0 | 0 | 0 | 1 | - | 0 | 0 | 0 | 0 | 0 | 0 | 0 | 0 |
| **h** | 0 | 1 | 0 | 0 | 0 | 0 | 0 | - | 0 | 0 | 0 | 0 | 0 | 0 | 0 |
| **i** | 0 | 0 | 0 | 0 | 0 | 0 | 0 | 1 | - | 0 | 1 | 0 | 0 | 0 | 0 |
| **j** | 0 | 0 | 0 | 0 | 0 | 0 | 0 | 0 | 1 | - | 0 | 0 | 0 | 0 | 0 |
| **k** | 0 | 0 | 0 | 0 | 0 | 0 | 0 | 0 | 0 | 1 | - | 0 | 0 | 0 | 0 |
| **l** | 0 | 0 | 1 | 0 | 0 | 0 | 0 | 0 | 0 | 0 | 0 | - | 0 | 0 | 0 |
| **m** | 0 | 0 | 0 | 0 | 0 | 0 | 0 | 0 | 0 | 0 | 0 | 1 | - | 0 | 1 |
| **n** | 0 | 0 | 0 | 0 | 0 | 0 | 0 | 0 | 0 | 0 | 0 | 0 | 1 | - | 0 |
| **o** | 0 | 0 | 0 | 0 | 0 | 0 | 0 | 0 | 0 | 0 | 0 | 0 | 0 | 1 | - |

According to the rule the alternatives *a*,*b*,*c* make a minimal weakly stable set. Thus, the alternatives *d*, *e*, *f*, *g*, *h*, *i*, *j*, *k*, *l*, *m*, *n*, *o* will be omitted on the first stage and the alternatives *a*,*b*,*c* will be chosen on the second stage of the two-stage choice procedure, i.e., $C(\vec{P}_X, X) = \{a, b, c\}$.

Suppose that the relative comparison of the alternative *a* and the alternative *d* is changed (*aμd*) while the relative comparison of any pair of other alternatives remains unchanged. Then a matrix of majority relation looks as

| | **a** | **b** | **c** | **d** | **e** | **f** | **g** | **h** | **i** | **j** | **k** | **l** | **m** | **n** | **o** |
|---|---|---|---|---|---|---|---|---|---|---|---|---|---|---|---|
| **a** | - | 0 | 1 | 0 | 0 | 0 | 1 | 1 | 0 | 0 | 0 | 0 | 0 | 0 | 0 |
| **b** | 1 | - | 0 | 0 | 0 | 0 | 0 | 0 | 0 | 0 | 1 | 1 | 0 | 0 | 0 |
| **c** | | 1 | - | 1 | 0 | 0 | 0 | 0 | 0 | 0 | 0 | 0 | 0 | 0 | 1 |
| **d** | 1 | 0 | 0 | - | 0 | 0 | 0 | 0 | 0 | 0 | 0 | 0 | 0 | 0 | 0 |
| **e** | 0 | 0 | 0 | 1 | - | 0 | 1 | 0 | 0 | 0 | 0 | 0 | 0 | 0 | 0 |
| **f** | 0 | 0 | 0 | 0 | 1 | - | 0 | 0 | 0 | 0 | 0 | 0 | 0 | 0 | 0 |
| **g** | 0 | 0 | 0 | 0 | 0 | 1 | - | 0 | 0 | 0 | 0 | 0 | 0 | 0 | 0 |
| **h** | 0 | 1 | 0 | 0 | 0 | 0 | 0 | - | 0 | 0 | 0 | 0 | 0 | 0 | 0 |
| **i** | 0 | 0 | 0 | 0 | 0 | 0 | 0 | 1 | - | 0 | 1 | 0 | 0 | 0 | 0 |
| **j** | 0 | 0 | 0 | 0 | 0 | 0 | 0 | 0 | 1 | - | 0 | 0 | 0 | 0 | 0 |
| **k** | 0 | 0 | 0 | 0 | 0 | 0 | 0 | 0 | 0 | 1 | - | 0 | 0 | 0 | 0 |
| **l** | 0 | 0 | 1 | 0 | 0 | 0 | 0 | 0 | 0 | 0 | 0 | - | 0 | 0 | 0 |
| **m** | 0 | 0 | 0 | 0 | 0 | 0 | 0 | 0 | 0 | 0 | 0 | 1 | - | 0 | 1 |
| **n** | 0 | 0 | 0 | 0 | 0 | 0 | 0 | 0 | 0 | 0 | 0 | 0 | 1 | - | 0 |
| **o** | 0 | 0 | 0 | 0 | 0 | 0 | 0 | 0 | 0 | 0 | 0 | 0 | 0 | 1 | - |

According to the rule the alternatives *a*,*b* make a minimal weakly stable set. Thus, the alternatives *c*, *d*, *e*, *f*, *g*, *h*, *i*, *j*, *k*, *l*, *m*, *n*, *o* will be omitted on the first stage and the alternative $b$ will be chosen on the second stage of the two-stage choice procedure (*bμa*), i.e., $C\left(\overrightarrow{P'_X}, X\right) = \{b\}$.

Then $\{a\} \in C(\vec{P}_X, X)$, $\{a\} \notin C\left(\overrightarrow{P'_X}, X\right)$. Thus, the Monotonicity conditionon 1 is not satisfied.

***383.2 Concordance condition (C)***

Let $X = \{a, b, c, d, e, f, g, h\}$ and a matrix of majority relation μ is the following

| | a | b | c | d | e | f | g | h |
|---|---|---|---|---|---|---|---|---|
| **a** | - | 1 | 0 | 1 | 0 | 0 | 0 | 0 |
| **b** | 0 | - | 1 | 1 | 1 | 0 | 0 | 0 |
| **c** | 1 | 0 | - | 0 | 0 | 0 | 0 | 1 |
| **d** | 0 | 0 | 1 | - | 0 | 0 | 0 | 0 |
| **e** | 0 | 0 | 1 | 0 | - | 0 | 0 | 0 |
| **f** | 0 | 0 | 0 | 0 | 1 | - | 0 | 0 |
| **g** | 0 | 0 | 0 | 0 | 0 | 1 | - | 0 |
| **h** | 0 | 0 | 0 | 0 | 0 | 0 | 1 | - |

On the first stage of the two-stage choice procedure minimal weakly stable set is defined. For this case there are 3 minimal weakly stable sets: $\{a, b\}, \{a, d\}$ and $\{b, c\}$. Thus, $Q = \{a, b\} \cup \{a, d\} \cup \{b, c\} = \{a, b, c, d\}$.

On the second stage of the two-stage choice procedure a Condorcet winner rule is used. For this case there is no Condorcet winner (*aμb*, *aμd*, *cμa*, *bμc*, *dμc*), Thus, $C\left(\vec{P}_X, X\right) = \emptyset$.

Consider now the subset $X' = X \backslash \{c, e, f, g\}$. A matrix of majority relation μ looks as

| | a | b | d |
|---|---|---|---|
| **a** | - | 1 | 1 |
| **b** | 0 | - | 1 |
| **d** | 0 | 0 | - |

On the first stage of the two-stage choice procedure minimal weakly stable set is defined. For this case the alternative *a* makes a minimal weakly stable set. Thus, $C\left(\vec{P}_{X'}, X'\right) = \{a\}$.

Finally, consider the subset $X'' = X \backslash \{b\}$. Then a matrix of majority relation μ looks as

| | a | c | d | e | f | g | h |
|---|---|---|---|---|---|---|---|
| **a** | - | 0 | 1 | 0 | 0 | 0 | 0 |
| **c** | 1 | - | 0 | 0 | 0 | 0 | 1 |
| **d** | 0 | 1 | - | 0 | 0 | 0 | 0 |
| **e** | 0 | 1 | 0 | - | 0 | 0 | 0 |
| **f** | 0 | 0 | 0 | 1 | - | 0 | 0 |
| **g** | 0 | 0 | 0 | 0 | 1 | - | 0 |
| **h** | 0 | 0 | 0 | 0 | 0 | 1 | - |

On the first stage of the two-stage choice procedure minimal weakly stable set is defined. For this case the alternatives *a* and *d* make a minimal weakly stable set. On the second stage of the two-stage choice procedure a Condorcet winner rule is used. Since *aμd*, $C\left(\vec{P}_{X''}, X''\right) = \{a\}$.

Then $C\left(\vec{P}_{X'}, X'\right) \cap C\left(\vec{P}_{X''}, X''\right) = \{a\} \nsubseteq C(\vec{P}_X, X)$. Thus, the condition **C** is not satisfied.

## Two-stage choice procedures 394-410, 412-420

Two-stage choice procedures do not satisfy the same conditions as choice procedures used on the second stage. To prove it we can use the same examples but with larger number of alternatives eliminated on the first stage of the two-stage choice procedure [11].

Thus, it is necessary to check those normative conditions of two-stage choice procedures which were satisfied for the second stage choice procedures.

- Two-stage choice procedure 412 does not satisfy the condition **H** (see paragraph 412.1 of Appendix 1).
- Two-stage choice procedures 404, 409 do not satisfy the condition **C**. To prove it the following example is used.

Let $X = \{a, b, c, d, e\}$ and a matrix of majority relation μ is the following

| | **a** | **b** | **c** | **d** | **e** | **f** |
|---|---|---|---|---|---|---|
| **a** | - | 1 | 0 | 0 | 1 | 1 |
| **b** | 0 | - | 1 | 1 | 0 | 1 |
| **c** | 1 | 0 | - | 1 | 1 | 0 |
| **d** | 1 | 0 | 0 | - | 0 | 0 |
| **e** | 0 | 0 | 0 | 1 | - | 1 |
| **f** | 0 | 0 | 1 | 1 | 0 | - |

Let us define the upper contour sets for each alternative. $D(a, \vec{P}_X) = \{c, d\}, D(b, \vec{P}_X) = \{a\}$, $D(c, \vec{P}_X) = \{b, f\}$, $D(d, \vec{P}_X) = \{b, c, e, f\}$, $D(e, \vec{P}_X) = \{a, c\}$, $D(f, \vec{P}_X) = \{a, b, e\}$. Then $c\gamma d$, $b\gamma e$, $b\gamma f$. Thus, $C(\vec{P}_X, X) = \{a, b, c\}$.

Consider now the subset $X' = X \backslash \{c, f\}$. Then a matrix of majority relation μ is the following

| | **a** | **b** | **d** | **e** |
|---|---|---|---|---|
| **a** | - | 1 | 0 | 1 |
| **b** | 0 | - | 1 | 0 |
| **d** | 1 | 0 | - | 0 |
| **e** | 0 | 0 | 1 | - |

Let us define the upper contour sets for each alternative. $D(a, \vec{P}_{X'}) = \{d\}, D(b, \vec{P}_{X'}) = \{a\}$, $D(d, \vec{P}_{X'}) = \{b, e\}$, $D(e, \vec{P}_{X'}) = \{a\}$. Thus, $C(\vec{P}_{X'}, X') = \{a, b, d, e\}$.

Finally, consider the subset $X'' = X \backslash \{a, b\}$. Then a matrix of majority relation μ is the following

| | **c** | **d** | **e** | **f** |
|---|---|---|---|---|
| **c** | - | 1 | 1 | 0 |
| **d** | 0 | - | 0 | 0 |
| **e** | 0 | 1 | - | 1 |
| **f** | 1 | 1 | 0 | - |

Let us define the upper contour sets for each alternative. $D(c, \vec{P}_{X''}) = \{f\}$, $D(d, \vec{P}_{X''}) = \{c, e, f\}$, $D(e, \vec{P}_{X''}) = \{c\}$, $D(f, \vec{P}_{X''}) = \{e\}$. Then $c\gamma d$, $e\gamma d$, $f\gamma d$. Thus, $C(\vec{P}_{X''}, X'') = \{c, e, f\}$.

Then $C(\vec{P}_{X'}, X') \cap C(\vec{P}_{X''}, X'') = \{e\} \nsubseteq C(\vec{P}_X, X)$. Thus, the condition **C** is not satisfied.

- Two-stage choice procedures 412 do not satisfy the condition **C** (see paragraph 412.2 of Appendix 1).
- Two-stage choice procedure 404 does not satisfy the condition **O**. To prove it the following example is used.

Let $X = \{a, b, c, d, e\}$ and a matrix of majority relation μ is the following

| | a | b | c | d | e |
|---|---|---|---|---|---|
| **a** | - | 1 | 0 | 1 | 1 |
| **b** | 0 | - | 1 | 1 | 0 |
| **c** | 1 | 0 | - | 1 | 0 |
| **d** | 0 | 0 | 0 | - | 1 |
| **e** | 0 | 0 | 0 | 0 | - |

Let us define the upper contour sets for each alternative. $D(a, \vec{P}_X) = \{c\}, D(b, \vec{P}_X) = \{a\}$, $D(c, \vec{P}_X) = \{b\}$, $D(d, \vec{P}_X) = \{a, b, c\}$, $D(e, \vec{P}_X) = \{a, d\}$. Then $c\gamma d$, $b\gamma e$. Thus, $C(\vec{P}_X, X) = \{a, b, c\}$.
Consider now the subset $X' = X \backslash \{d\}$. Then a matrix of majority relation μ is the following

| | a | b | c | e |
|---|---|---|---|---|
| **a** | - | 1 | 0 | 1 |
| **b** | 0 | - | 1 | 0 |
| **c** | 1 | 0 | - | 0 |
| **e** | 0 | 0 | 0 | - |

Let us define the upper contour sets for each alternative. $D(a, \vec{P}_{X'}) = \{c\}, D(b, \vec{P}_{X'}) = \{a\}$, $D(c, \vec{P}_{X'}) = \{b\}$, $D(e, \vec{P}_{X'}) = \{a\}$. Thus, $C(\vec{P}_{X'}, X') = \{a, b, c, e\}$.
Then $C(\vec{P}_{X'}, X') \neq C(\vec{P}_X, X)$, Consequently, the condition **O** is not satisfied.

- Two-stage choice procedures 394-403, 405-407, 409, 412, 414-420 do not satisfy the Monotonicity conditionon 1 (the proof follows from the paragraphs 394.5, 409.5 of Appendix 1 and properties of second stage choice procedures).
- Two-stage choice procedure 412 does not satisfy the Monotonicity conditionon 2 (see paragraph 412.1 of Appendix 1).
- Two-stage choice procedures 393-420 do not satisfy the non-compensatory condition (see the properties of Fishburn rule).

***394.5 Monotonicity condition 1.***
Consider the example from paragraph 412.1 of Appendix 1. Suppose that $b\mu c$. Then the alternative $f$ is eliminated on the first stage and, consequently, the alternative $b$ is not chosen on the second stage. Since $\{b\} \in C(\vec{P}_X, X)$, $\{b\} \notin C\left(\overrightarrow{P'_X}, X\right)$, the Monotonicity conditionon 1 is not satisfied.

***409.5 Monotonicity condition 1.***
Let $X = \{a, b, c, d, e, f\}$ and a matrix of majority relation μ is the following

| | a | b | c | d | e | f |
|---|---|---|---|---|---|---|
| **a** | - | 0 | 1 | 0 | 1 | 0 |
| **b** | 0 | - | 0 | 1 | 0 | 1 |
| **c** | 0 | 0 | - | 0 | 0 | 0 |
| **d** | 0 | 0 | 0 | - | 1 | 0 |
| **e** | 0 | 1 | 0 | 0 | - | 1 |
| **f** | 1 | 0 | 0 | 0 | 0 | - |

According to the rule the alternatives $a$,$b$ will be chosen, i.e., $C(\vec{P}_X, X) = \{a, b\}$.

Suppose that the relative position of the alternative $a$ is improved such that $a\mu b$, $a\mu d$ while the relative comparison of any pair of other alternatives remains unchanged. Then a matrix of majority relation looks as

| | **a** | **b** | **c** | **d** | **e** | **f** |
|---|---|---|---|---|---|---|
| **a** | - | 1 | 1 | 1 | 1 | 0 |
| **b** | 0 | - | 0 | 1 | 0 | 1 |
| **c** | 0 | 0 | - | 0 | 0 | 0 |
| **d** | 0 | 0 | 0 | - | 1 | 0 |
| **e** | 0 | 1 | 0 | 0 | - | 1 |
| **f** | 1 | 0 | 0 | 0 | 0 | - |

For this case the alternatives $a$,$c$,$f$ are considered for the second stage of the choice procedure. Since $f\mu a$, $a\mu c$, $C\left(\overrightarrow{P'_X}, X\right) = \{f, c\}$.

Then $\{a\} \in C(\vec{P}_X, X)$, $\{a\} \notin C\left(\overrightarrow{P'_X}, X\right)$. Thus, the Monotonicity conditionon 1 is not satisfied.

***412.1 Heredity condition (H)***

Let $X = \{a, b, c, d, e, f\}$ and a matrix of majority relation μ is the following

| | **a** | **b** | **c** | **d** | **e** | **f** |
|---|---|---|---|---|---|---|
| **a** | - | 1 | 1 | 0 | 1 | 0 |
| **b** | 0 | - | 0 | 1 | 0 | 1 |
| **c** | 0 | 0 | - | 0 | 0 | 1 |
| **d** | 0 | 0 | 0 | - | 1 | 0 |
| **e** | 0 | 1 | 0 | 0 | - | 0 |
| **f** | 1 | 0 | 0 | 0 | 0 | - |

According to the two-stage choice procedure the alternatives $a$,$b$ will be chosen, i.e., $C(\vec{P}_X, X) = \{a, b\}$.

Consider now the subset $X' = X \backslash \{c, d\}$. A matrix of majority relation μ looks as

| | **a** | **b** | **e** | **f** |
|---|---|---|---|---|
| **a** | - | 0 | 1 | 0 |
| **b** | 0 | - | 0 | 1 |
| **e** | 0 | 1 | - | 0 |
| **f** | 1 | 0 | 0 | - |

According to the rule the choice is empty, i.e., $C(\vec{P}_{X'}, X') = \emptyset$.

Then $C(\vec{P}_X, X) \cap X' = \{a, b\} \not\subseteq C(\vec{P}_{X'}, X')$. Thus, the condition **H** is not satisfied.

*412.2 Concordance condition (C)*

Let $X = \{a, b, c, d, e, f, g\}$ and a matrix of majority relation μ is the following

| | a | b | c | d | e | f | g |
|---|---|---|---|---|---|---|---|
| **a** | - | 0 | 1 | 0 | 1 | 0 | 0 |
| **b** | 0 | - | 0 | 1 | 0 | 1 | 1 |
| **c** | 0 | 0 | - | 0 | 0 | 1 | 0 |
| **d** | 0 | 0 | 0 | - | 1 | 0 | 0 |
| **e** | 0 | 1 | 0 | 0 | - | 0 | 1 |
| **f** | 1 | 0 | 0 | 0 | 0 | - | 0 |
| **g** | 1 | 0 | 0 | 1 | 0 | 1 | - |

According to the two-stage choice procedure the alternatives *a,b* will be chosen, i.e., $C(\vec{P}_X, X) = \{a, b\}$.

Consider now the subset $X' = X \backslash \{b\}$. A matrix of majority relation μ looks as

| | a | c | d | e | f | g |
|---|---|---|---|---|---|---|
| **a** | - | 1 | 0 | 1 | 0 | 0 |
| **c** | 0 | - | 0 | 0 | 1 | 0 |
| **d** | 0 | 0 | - | 1 | 0 | 0 |
| **e** | 0 | 0 | 0 | - | 0 | 1 |
| **f** | 1 | 0 | 0 | 0 | - | 0 |
| **g** | 1 | 0 | 1 | 0 | 1 | - |

According to the given two-stage choice procedure the alternatives *c,g* will be chosen, i.e., $C(\vec{P}_{X'}, X') = \{c, g\}$.

Finally, consider the subset $X'' = X \backslash \{a\}$. Then a matrix of majority relation μ looks as

| | b | c | d | e | f | g |
|---|---|---|---|---|---|---|
| **b** | - | 0 | 1 | 0 | 1 | 1 |
| **c** | 0 | - | 0 | 0 | 1 | 0 |
| **d** | 0 | 0 | - | 1 | 0 | 0 |
| **e** | 1 | 0 | 0 | - | 0 | 1 |
| **f** | 0 | 0 | 0 | 0 | - | 0 |
| **g** | 0 | 0 | 1 | 0 | 1 | - |

According to the rule the alternative *c* will be chosen, i.e., $C(\vec{P}_{X''}, X'') = \{c\}$.

Then $C(\vec{P}_{X'}, X') \cap C(\vec{P}_{X''}, X'') = \{c\} \not\subseteq C(\vec{P}_X, X)$. Thus, the condition **C** is not satisfied.

*422.5 Monotonicity condition 1.*

Let $X = \{a, b, c, d, e, f\}$ and a matrix of majority relation μ is the following

| | a | b | c | d | e | f |
|---|---|---|---|---|---|---|
| **a** | - | 1 | 0 | 0 | 0 | 0 |
| **b** | 0 | - | 1 | 0 | 0 | 0 |
| **c** | 1 | 0 | - | 1 | 1 | 0 |
| **d** | 0 | 0 | 0 | - | 0 | 0 |
| **e** | 0 | 0 | 0 | 0 | - | 0 |
| **f** | 1 | 0 | 0 | 1 | 0 | - |

According to the rule the alternatives *a,b,c* will be chosen, i.e., $C(\vec{P}_X, X) = \{a, b, c\}$.

Suppose that the relative position of the alternative $a$ is improved such that $a\mu c$ while the relative comparison of any pair of other alternatives remains unchanged. Then a matrix of majority relation looks as

| | **a** | **b** | **c** | **d** | **e** | **f** |
|---|---|---|---|---|---|---|
| **a** | - | 1 | 0 | 0 | 0 | 0 |
| **b** | 0 | - | 1 | 0 | 0 | 0 |
| **c** | 0 | 0 | - | 1 | 1 | 0 |
| **d** | 0 | 0 | 0 | - | 0 | 0 |
| **e** | 0 | 0 | 0 | 0 | - | 0 |
| **f** | 1 | 0 | 0 | 1 | 0 | - |

For this case the alternatives *a,b,c,f* are considered for the second stage of the choice procedure. For this majority matrix it is possible to construct a profile $\overrightarrow{P'_X}$ such that $\{a\} \notin C\left(\overrightarrow{P'_X}, X\right)$.

Then $\{a\} \in C(\vec{P}_X, X)$, $\{a\} \notin C\left(\overrightarrow{P'_X}, X\right)$. Thus, the Monotonicity conditionon 1 is not satisfied.

***440.1 Heredity condition (H)***

Let $X = \{a, b, c, d, e, f\}$ and a matrix of majority relation μ is the following

| | **a** | **b** | **c** | **d** | **e** | **f** |
|---|---|---|---|---|---|---|
| **a** | - | 0 | 1 | 0 | 0 | 0 |
| **b** | 0 | - | 1 | 0 | 0 | 0 |
| **c** | 0 | 0 | - | 1 | 1 | 1 |
| **d** | 0 | 0 | 0 | - | 1 | 0 |
| **e** | 0 | 0 | 0 | 0 | - | 1 |
| **f** | 0 | 0 | 0 | 1 | 0 | - |

According to the two-stage choice procedure the alternatives *a,b* will be chosen, i.e., $C(\vec{P}_X, X) = \{a, b\}$.

Consider now the subset $X' = X \backslash \{c\}$. A matrix of majority relation μ looks as

| | **a** | **b** | **d** | **e** | **f** |
|---|---|---|---|---|---|
| **a** | - | 0 | 0 | 0 | 0 |
| **b** | 0 | - | 0 | 0 | 0 |
| **d** | 0 | 0 | - | 1 | 0 |
| **e** | 0 | 0 | 0 | - | 1 |
| **f** | 0 | 0 | 1 | 0 | - |

According to the rule the choice is empty, i.e., $C(\vec{P}_{X'}, X') = \emptyset$.

Then $C(\vec{P}_X, X) \cap X' = \{a, b\} \not\subseteq C(\vec{P}_{X'}, X')$. Thus, the condition **H** is not satisfied.

***465.2 Concordance condition (C)***

Consider the example from paragraph 468.2 of Appendix 1. Then $C(\vec{P}_X, X) = \{a, b, c, g\}$.

Consider now the subset $X' = X\setminus\{b\}$. A matrix of majority relation μ looks as

| | **a** | **b** | **c** | **d** | **e** | **f** |
|---|---|---|---|---|---|---|
| **a** | - | 0 | 0 | 0 | 1 | 1 |
| **b** | 0 | - | 0 | 0 | 0 | 0 |
| **c** | 1 | 0 | - | 0 | 0 | 1 |
| **d** | 0 | 0 | 1 | - | 0 | 1 |
| **e** | 0 | 1 | 1 | 1 | - | 0 |
| **f** | 0 | 0 | 0 | 0 | 1 | - |

According to the given two-stage choice procedure the alternatives *a*,*b*,*d* will be chosen, i.e., $C(\vec{P}_{X'}, X') = \{a, b, d\}$.

Finally, consider the subset $X'' = X\setminus\{a\}$. Then a matrix of majority relation μ looks as

| | **b** | **c** | **d** | **e** | **f** | **g** |
|---|---|---|---|---|---|---|
| **b** | - | 0 | 0 | 0 | 0 | 0 |
| **c** | 0 | - | 0 | 0 | 1 | 0 |
| **d** | 0 | 1 | - | 0 | 1 | 0 |
| **e** | 1 | 1 | 1 | - | 0 | 0 |
| **f** | 0 | 0 | 0 | 1 | - | 1 |
| **g** | 0 | 0 | 1 | 1 | 0 | - |

According to the rule the alternatives *b*,*c*,*d*,*g* will be chosen, i.e., $C(\vec{P}_{X''}, X'') = \{b, c, d, g\}$.

Then $C(\vec{P}_{X'}, X') \cap C(\vec{P}_{X''}, X'') = \{b, d\} \not\subseteq C(\vec{P}_X, X)$. Thus, the condition **C** is not satisfied.

***465.5 Monotonicity condition 1.***

Let $X = \{a, b, c, d\}$ and a matrix of majority relation μ is the following

| | **a** | **b** | **c** | **d** |
|---|---|---|---|---|
| **a** | - | 1 | 0 | 0 |
| **b** | 0 | - | 0 | 1 |
| **c** | 1 | 0 | - | 0 |
| **d** | 0 | 0 | 1 | - |

According to the rule the alternatives *a*,*b*,*c*,*d* will be chosen. Suppose that the relative position of the alternative $a$ is improved ($a\mu d$). Then the alternative *d* is omitted on the first stage. For this case it is possible to construct a profile $\overrightarrow{P'_X}$ such that $\{a\} \notin C\left(\overrightarrow{P'_X}, X\right)$.

Then $\{a\} \in C(\vec{P}_X, X)$, $\{a\} \notin C\left(\overrightarrow{P'_X}, X\right)$. Thus, the Monotonicity conditionon 1 is not satisfied.

### *468.1 Heredity condition (H)*

Let $X = \{a, b, c, d, e, f, g\}$ and a matrix of majority relation μ is the following

| | a | b | c | d | e | f | g |
|---|---|---|---|---|---|---|---|
| **a** | - | 0 | 0 | 0 | 1 | 1 | 0 |
| **b** | 0 | - | 0 | 0 | 0 | 0 | 0 |
| **c** | 1 | 0 | - | 0 | 0 | 1 | 0 |
| **d** | 0 | 0 | 1 | - | 0 | 1 | 0 |
| **e** | 0 | 1 | 1 | 1 | - | 0 | 0 |
| **f** | 0 | 0 | 0 | 0 | 1 | - | 1 |
| **g** | 0 | 0 | 0 | 0 | 1 | 0 | - |

According to the given two-stage choice procedure the alternatives *a,b,d,g* will be chosen, i.e., $C(\vec{P}_X, X) = \{a, b, d, g\}$.

Consider now the subset $X' = X\backslash\{d, e\}$. A matrix of majority relation μ looks as

| | a | b | c | f | g |
|---|---|---|---|---|---|
| **a** | - | 0 | 0 | 1 | 0 |
| **b** | 0 | - | 0 | 0 | 0 |
| **c** | 1 | 0 | - | 1 | 0 |
| **f** | 0 | 0 | 0 | - | 1 |
| **g** | 0 | 0 | 0 | 0 | - |

According to the rule the alternatives *b,c,g* will be chosen, i.e., $C(\vec{P}_{X'}, X') = \{b, c, g\}$.

Then $C(\vec{P}_X, X) \cap X' = \{a, b, g\} \not\subseteq C(\vec{P}_{X'}, X')$. Thus, the condition **H** is not satisfied.

### *468.2 Concordance condition (C)*

Let $X = \{a, b, c, d, e, f, g\}$ and a matrix of majority relation μ is the following

| | a | b | c | d | e | f | g |
|---|---|---|---|---|---|---|---|
| **a** | - | 0 | 0 | 0 | 1 | 1 | 0 |
| **b** | 0 | - | 0 | 0 | 0 | 0 | 0 |
| **c** | 1 | 0 | - | 0 | 0 | 1 | 0 |
| **d** | 0 | 0 | 1 | - | 0 | 1 | 0 |
| **e** | 0 | 1 | 1 | 1 | - | 0 | 0 |
| **f** | 0 | 0 | 0 | 0 | 1 | - | 1 |
| **g** | 0 | 0 | 0 | 1 | 1 | 0 | - |

According to the two-stage choice procedure the alternatives *b,g* will be chosen, i.e., $C(\vec{P}_X, X) = \{b, g\}$.

Consider now the subset $X' = X\backslash\{d\}$. A matrix of majority relation μ looks as

| | a | b | c | e | f | g |
|---|---|---|---|---|---|---|
| **a** | - | 0 | 0 | 1 | 1 | 0 |
| **b** | 0 | - | 0 | 0 | 0 | 0 |
| **c** | 1 | 0 | - | 0 | 1 | 0 |
| **e** | 0 | 1 | 1 | - | 0 | 0 |
| **f** | 0 | 0 | 0 | 1 | - | 1 |
| **g** | 0 | 0 | 0 | 1 | 0 | - |

According to the given two-stage choice procedure the alternatives *b,c,g* will be chosen, i.e., $C(\vec{P}_{X'}, X') = \{b, c, g\}$.

Finally, consider the subset $X'' = X\backslash\{a,b,e,f\}$. Then a matrix of majority relation μ looks as

| | **c** | **d** | **g** |
|---|---|---|---|
| **c** | - | 0 | 0 |
| **d** | 1 | - | 0 |
| **g** | 0 | 1 | - |

According to the rule the alternatives *c*,*g* will be chosen, i.e., $C(\vec{P}_{X''}, X'') = \{c,g\}$.

Then $C(\vec{P}_{X'}, X') \cap C(\vec{P}_{X''}, X'') = \{c,g\} \not\subseteq C(\vec{P}_X, X)$. Thus, the condition **C** is not satisfied.

***468.5 Monotonicity condition 1.***

Let $X = \{a,b,c,d,e,f,g\}$ and a matrix of majority relation μ is the following

| | **a** | **c** | **d** | **e** | **f** | **g** |
|---|---|---|---|---|---|---|
| **a** | - | 0 | 0 | 0 | 1 | 0 |
| **c** | 1 | - | 0 | 0 | 1 | 0 |
| **d** | 0 | 1 | - | 0 | 1 | 0 |
| **e** | 0 | 1 | 1 | - | 0 | 0 |
| **f** | 0 | 0 | 0 | 1 | - | 1 |
| **g** | 0 | 0 | 1 | 1 | 0 | - |

According to the rule the alternatives *e*,*f* are omitted on the first stage. Suppose that the relative position of the alternative $a$ is improved such that $a\mu g$ while the relative comparison of any pair of other alternatives remains unchanged. Then the alternatives *e*,*f* are omitted on the first stage. For this case it is possible to construct a profile $\overrightarrow{P_X'}$ such that $\{a\} \notin C\left(\overrightarrow{P_X'}, X\right)$.

Then $\{a\} \in C(\vec{P}_X, X)$, $\{a\} \notin C\left(\overrightarrow{P_X'}, X\right)$. Thus, the Monotonicity conditionon 1 is not satisfied.

## Two-stage choice procedures 478-494, 496-504

Two-stage choice procedures 478-494, 496-504 satisfy the same conditions as two-stage choice procedures which use the Fishburn rule and Uncovered set I on the first stage. To prove it we can use the same examples that are used for these rules [11].

## Two-stage choice procedures 534-538,543,554,558

Two-stage choice procedures do not satisfy the same conditions as choice procedures used on the second stage. As for other normative conditions, two-stage choice procedures 534-538,543,554,558 do not satisfy any of them as the subset of alternatives for the second stage can be easily changed when the core is applied on the first stage of the rule. To prove it we can provide the following example.

***534.5 Monotonicity condition 1.***

Let $X = \{a,b,c,d\}$ and the profile $\vec{P}_X$ is the following

| $P_1$ | $P_2$ | $P_3$ | $P_4$ | $P_5$ | $P_6$ | $P_7$ | $P_8$ |
|---|---|---|---|---|---|---|---|
| *a* | *a* | *a* | *c* | *b* | *b* | *d* | *d* |
| *c* | *c* | *c* | *a* | *d* | *d* | *b* | *b* |
| *b* | *d* | *d* | *b* | *c* | *a* | *c* | *c* |
| *d* | *b* | *b* | *d* | *a* | *c* | *a* | *a* |

According to the rule the alternative $a$ will be chosen, i.e., $C(\vec{P}_X, X) = \{a\}$.

Suppose now that the position of the alternative $a$ was improved in $P'_6$ while the relative comparison of any pair of other alternatives remained unchanged.

| $P'_1$ | $P'_2$ | $P'_3$ | $P'_4$ | $P'_5$ | $P'_6$ | $P'_7$ | $P'_8$ |
|---|---|---|---|---|---|---|---|
| $a$ | $a$ | $a$ | $c$ | $b$ | $b$ | $d$ | $d$ |
| $c$ | $c$ | $c$ | $a$ | $d$ | $a$ | $b$ | $b$ |
| $b$ | $d$ | $d$ | $b$ | $c$ | $d$ | $c$ | $c$ |
| $d$ | $b$ | $b$ | $d$ | $a$ | $c$ | $a$ | $a$ |

Then the alternative $d$ will be omitted after applying the first choice procedure (core) and the subset $X' = X\backslash\{d\}$ will be presented for the second choice procedure (the plurality rule). A contraction of a profile $\overrightarrow{P'_X}$ onto a set $X'$, i.e., $\overrightarrow{P'}_{X'}$, looks as

| $P'_1$ | $P'_2$ | $P'_3$ | $P'_4$ | $P'_5$ | $P'_6$ | $P'_7$ | $P'_8$ |
|---|---|---|---|---|---|---|---|
| $a$ | $a$ | $a$ | $c$ | $b$ | $b$ | $b$ | $b$ |
| $c$ | $c$ | $c$ | $a$ | $c$ | $a$ | $c$ | $c$ |
| $b$ | $b$ | $b$ | $b$ | $a$ | $c$ | $a$ | $a$ |

According to the rule the alternative $b$ will be chosen, i.e., $C\left(\overrightarrow{P'_X}, X\right) = C\left(\overrightarrow{P'_{X'}}, X'\right) = \{b\}$.

Then $\{a\} \in C(\vec{P}_X, X)$, $\{a\} \notin C\left(\overrightarrow{P'_X}, X\right)$. Thus, the Monotonicity conditionon 1 is not satisfied for the given two-stage choice procedure.

## Two-stage choice procedures 561-588

Two-stage choice procedures 561-588 do not satisfy the same normative conditions as two-stage choice procedures which use minimal weakly stable set on the first stage. To prove it similar examples can be used but with larger number of alternatives (see the properties of k-stable set) [11].

## Two-stage choice procedures 589-616

To check the properties of two-stage choice procedures 589-616 it is possible to use the same counter-examples that were used to check the properties of the threshold rule [11]. Thus, it remained to check normative conditions which are satisfied for the threshold rule.

- Two-stage choice procedures 589-616 satisfy the Monotonicity conditionon 1. The proof follows from properties of the threshold rule (an improved position of any chosen alternative $x$ leads to the choice of only this alternative on the first stage of two-stage procedure).
- Two-stage choice procedures 589-616 do not satisfy the non-compensatory condition. To proof it the following example is provided.

Let $X = \{a, b, c, d\}$ and the profile $\vec{P}_X$ is the following

| $P_1$ | $P_2$ | $P_3$ |
|---|---|---|
| $a$ | $b$ | $c$ |
| $d$ | $a$ | $b$ |
| $b$ | $c$ | $a$ |
| $c$ | $d$ | $d$ |

According to the first stage rule, the alternatives $a$ and $b$ will be chosen on the first stage of the two-stage choice procedure. Consider now the subset $X' = X\backslash\{c, d\}$. A contraction of a profile $\vec{P}_X$ onto a set $X'$, i.e., $\vec{P}_{X'}$, looks as

| **P₁** | **P₂** | **P₃** |
|---|---|---|
| $a$ | $b$ | $b$ |
| $b$ | $a$ | $a$ |

According to the second stage rule, the alternative $b$ will be chosen on the second stage of the two-stage choice procedure. Thus, $C(\vec{P}_X, X) = \{b\}$.

Now let us write the profile $\vec{P}_X$ in the following form

| $X$ | $\varphi_1$ | $\varphi_2$ | $\varphi_3$ |
|---|---|---|---|
| $a$ | 4 | 3 | 2 |
| $b$ | 2 | 4 | 3 |
| $c$ | 1 | 2 | 4 |
| $d$ | 3 | 1 | 1 |

According to the non-compensatory condition the alternatives $a$,$b$ (the alternatives $a$ and $b$ are equal) are better than the alternative $c$ and the alternative $c$ is better than the alternative $d$. Thus, the non-compensatory condition is not satisfied as $\{a, b\} \neq C(\vec{P}_X, X)$.

## Two-stage choice procedures 617-700

Two-stage choice procedures do not satisfy the same conditions as Copeland rules 1-3 and choice procedures used on the second stage. To prove it similar examples can be used but with larger number of alternatives [11]. Thus, it remained to check normative conditions which are satisfied for Copeland rules 1-3 and choice procedures used on the second stage.

- Two-stage choice procedures 617-620, 623-624, 628-648, 651-652, 656-676, 679-680, 684-700 satisfy the Monotonicity conditionon 1 (the proof follows from properties of Copeland rules 1-3 and choice procedures used on the second stage).
- Two-stage choice procedures 621-622, 625-627, 649-650, 653-655, 677-678, 681-683 do not satisfy the Monotonicity conditionon 1 (the same example that was used to check the properties of the second stage choice procedures can be provided).
- Two-stage choice procedures 617-700 do not satisfy the non-compensatory condition (the proof follows from properties of Copeland rules 1-3).

## Two-stage choice procedures 701-728

Two-stage choice procedures do not satisfy the same conditions as choice procedures used on the second stage. To prove it similar examples can be used but with larger number of alternatives [11]. Thus, it remained to check normative conditions which areaz satisfied for choice procedures used on the second stage.

- Two-stage choice procedures 701, 719-720 do not satisfy the condition **H** (the proof follows from properties of super-threshold rule).
- Two-stage choice procedures 712, 717, 719-720 do not satisfy the condition **C** (the proof follows from properties of super-threshold rule).
- Two-stage choice procedures 712, 717 do not satisfy the condition **O** (the proof follows from properties of super-threshold rule).

- Two-stage choice procedures 701, 712-713, 719-720 satisfy the Monotonicity conditionon 1 (the proof follows from properties of threshold-rule and choice procedures used on the second stage).
- Two-stage choice procedures 702-711, 714-718, 721-728 do not satisfy the Monotonicity conditionon 1. Since super-threshold rule narrows the subset of alternatives remained after the first stage of the procedure, alternatives dominated by the chosen one might be eliminated on the first stage and, consequently, on the second stage of the rule the chosen alternative will not be included into a new choice.
- Two-stage choice procedures 701-728 do not satisfy the non-compensatory condition (the proof follows from properties of super-threshold rule).

## Two-stage choice procedures 729-784

To check the properties of two-stage choice procedures which use minimax and Simpson procedures on the first stage we can use the same counter-examples which were used to check the properties of minimax and Simpson procedures [11]. Thus, it remained to study those normative conditions which are satisfied for choice procedures used on the first stage.

- Two-stage choice procedures 729, 740-741, 747-748, 757, 768-769, 775-776 satisfy the Monotonicity conditionon 1 (the proof follows from properties of Minimax and Simpson procedures and choice procedures used on the second stage).
- Two-stage choice procedures 730-739, 742-746, 749-756, 758-767, 770-774, 777-784 do not satisfy the Monotonicity conditionon 1. To prove it the following example for the two-stage choice procedure 731 is given.

Let $X = \{a, b, c\}$ and the profile $\vec{P}_X$ looks as

| $P_1$ | $P_2$ | $P_3$ | $P_4$ | $P_5$ | $P_6$ | $P_7$ | $P_8$ | $P_9$ | $P_{10}$ | $P_{11}$ | $P_{12}$ | $P_{13}$ | $P_{14}$ | $P_{15}$ |
|---|---|---|---|---|---|---|---|---|---|---|---|---|---|---|
| *c* | *c* | *b* | *a* | *a* | *a* | *a* | *a* | *c* | *d* | *d* | *d* | *d* | *b* | *c* |
| *b* | *b* | *c* | *b* | *c* | *c* | *c* | *c* | *b* | *c* | *b* | *b* | *b* | *a* | *d* |
| *d* | *d* | *d* | *c* | *b* | *b* | *b* | *b* | *d* | *a* | *a* | *a* | *a* | *d* | *a* |
| *a* | *a* | *a* | *d* | *d* | *d* | *d* | *d* | *a* | *b* | *c* | *c* | *c* | *c* | *b* |

Let us construct a matrix $S^{-}(\vec{P}_X, X)$ for the profile $\vec{P}_X$.

| | **a** | **b** | **c** | **d** |
|---|---|---|---|---|
| **a** | - | 7 | 9 | 6 |
| **b** | 8 | - | 6 | 10 |
| **c** | 6 | 9 | - | 10 |
| **d** | 9 | 5 | 5 | - |

According to the minimax procedure used on the first stage, the alternatives *a,b,c* will be chosen, i.e., $C_1(\vec{P}_X, X) = \{a, b, c\}$, where $C_1$ is a Minimax choice procedure.

Consider now the subset $X' = X\backslash\{d\}$. A contraction of a profile $\vec{P}_X$ onto a set $X'$, i.e., $\vec{P}_{X'}$, looks as

| $P_1$ | $P_2$ | $P_3$ | $P_4$ | $P_5$ | $P_6$ | $P_7$ | $P_8$ | $P_9$ | $P_{10}$ | $P_{11}$ | $P_{12}$ | $P_{13}$ | $P_{14}$ | $P_{15}$ |
|---|---|---|---|---|---|---|---|---|---|---|---|---|---|---|
| *c* | *c* | *b* | *a* | *a* | *a* | *a* | *a* | *c* | *c* | *b* | *b* | *b* | *b* | *c* |
| *b* | *b* | *c* | *b* | *c* | *c* | *c* | *c* | *b* | *a* | *a* | *a* | *a* | *a* | *a* |
| *a* | *a* | *a* | *c* | *b* | *b* | *b* | *b* | *a* | *b* | *c* | *c* | *c* | *c* | *b* |

According to the rule the alternative *a* will be chosen, i.e., $C(\vec{P}_X, X) = \{a\}$.

Consider now a profile $\overrightarrow{P'_X}$, which differs from the profile $\vec{P}_X$ only by improved position of the alternative $a$ in $P'_{10}$.

| $P'_1$ | $P'_2$ | $P'_3$ | $P'_4$ | $P'_5$ | $P'_6$ | $P'_7$ | $P'_8$ | $P'_9$ | $P'_{10}$ | $P'_{11}$ | $P'_{12}$ | $P'_{13}$ | $P'_{14}$ | $P'_{15}$ |
|---|---|---|---|---|---|---|---|---|---|---|---|---|---|---|
| *c* | *c* | *b* | *a* | *a* | *a* | *a* | *a* | *c* | *d* | *d* | *d* | *d* | *b* | *c* |
| *b* | *b* | *c* | *b* | *c* | *c* | *c* | *c* | *b* | *a* | *b* | *b* | *b* | *a* | *d* |
| *d* | *d* | *d* | *c* | *b* | *b* | *b* | *b* | *d* | *c* | *a* | *a* | *a* | *d* | *a* |
| *a* | *a* | *a* | *d* | *d* | *d* | *d* | *d* | *a* | *b* | *c* | *c* | *c* | *c* | *b* |

Let us construct a matrix $S^{-}(\overrightarrow{P'}_X, X)$ for the profile $\overrightarrow{P'}_X$.

| | **a** | **b** | **c** | **d** |
|---|---|---|---|---|
| **a** | - | 7 | 10 | 6 |
| **b** | 8 | - | 6 | 10 |
| **c** | 5 | 9 | - | 10 |
| **d** | 9 | 5 | 5 | - |

According to the minimax procedure used on the first stage, the alternatives *a,b* will be chosen, i.e., $C_1\left(\overrightarrow{P'}_X, X\right) = \{a, b\}$, where $C_1$ is a Minimax choice procedure.

Consider now the subset $X'' = X\backslash\{c, d\}$. A contraction of a profile $\overrightarrow{P'}_X$ onto a set $X''$, i.e., $\overrightarrow{P'}_{X''}$, looks as

| $P'_1$ | $P'_2$ | $P'_3$ | $P'_4$ | $P'_5$ | $P'_6$ | $P'_7$ | $P'_8$ | $P'_9$ | $P'_{10}$ | $P'_{11}$ | $P'_{12}$ | $P'_{13}$ | $P'_{14}$ | $P'_{15}$ |
|---|---|---|---|---|---|---|---|---|---|---|---|---|---|---|
| *b* | *b* | *b* | *a* | *a* | *a* | *a* | *a* | *b* | *a* | *b* | *b* | *b* | *b* | *a* |
| *a* | *a* | *a* | *b* | *b* | *b* | *b* | *b* | *a* | *b* | *a* | *a* | *a* | *a* | *b* |

According to the rule the alternative *b* will be chosen, i.e., $C\left(\overrightarrow{P'}_X, X\right) = \{b\}$.

Then $\{a\} \in C(\vec{P}_X, X)$, $\{a\} \notin C\left(\overrightarrow{P'_X}, X\right)$. Thus, the Monotonicity conditionon 1 is not satisfied.

## References

1. Aizerman M., Aleskerov F. Theory of Choice. Elsevier, North-Holland, 1995.
2. Aleskerov F. Procedures of multicriterial choice. Preprints of the IFAC/IFORS Conference on Control Science and Technology for Development, Beijing, China, 1985. P. 858–869.
3. Aleskerov F. Multicriterial interval choice models // Information Sciences. 1994. No. 1. P. 14–26.
4. Aleskerov F.T., Khabina E.L., Shvarts D.A. Binary relations, graphs and collective decisions (in Russian). Moscow: HSE Publishing House, 2006.
5. Aleskerov F.T., Kurbanov E. On the degree of manipulability of group choice rules // Automation and Remote Control (in Russian). 1998. No. 10. P. 134–146.
6. Volsky V.I. Voting rules in small groups from ancient times to the XX century (in Russian). Working paper WP7/2014/02. Moscow: HSE Publishing House, 2014.
7. Aleskerov F.T., Subochev A. Matrix-vector representation of various solution concepts. Working papers WP7/2009/03. Moscow: HSE Publishing House, 2009.
8. Subochev A. Dominant, Weakly Stable, Uncovered Sets: Properties and Extensions. Working papers WP7/2008/03. Moscow: HSE Publishing House, 2008.
9. Aleskerov F.T., Yuzbashev D.V., Yakuba V.I. Threshold aggregation for three-graded rankings // Automation and Remote Control. 2007. No. 1. P. 147–152.
10. Aleskerov F., Chistyakov V., Kalyagin V. Social threshold aggregations // Social Choice and Welfare. 2010. Vol. 35. No. 4. P. 627–646.
11. Shvydun S. Normative properties of multi-criteria choice procedures and their superpositions: I. Working paper WP7/2015/07 (Part 1). Moscow: HSE Publishing House, 2015.
12. Aleskerov F., Cinar Y. 'q-Pareto-scalar' Two-stage Exstremization Model and its Reducibility to One-stage Model // Theory and Decision. 2008. No. 65. P. 291–304.
13. Aleskerov F. Arrovian Aggregation Models, Kluwer Academic Publishers, Dordrecht, Boston. L., 1999.

Исследуются двухступенчатые процедуры выбора, которые представляют собой суперпозицию двух процедур выбора. Показано, какие из рассматриваемых процедур выбора удовлетворяют существующим нормативным условиям, описывающим, каким образом изменяется конечный выбор при изменении предъявляемого множества альтернатив и оценок альтернатив по критериям. Особое внимание уделяется двухступенчатым процедурам, в основе которых лежат позиционные правила, а также правила, использующие мажоритарное отношение, вспомогательную числовую шкалу и турнирную матрицу. Приводится теорема о том, какие нормативные условия выполняются для рассматриваемых двухступенчатых процедур выбора. Оценена вычислительная сложность двухступенчатых процедур выбора и время их выполнения на реальных данных.

*Швыдун С.В.* – НИУ ВШЭ, ИПУ РАН.

Швыдун С. В.

**Нормативные свойства процедур многокритериального выбора и их суперпозиции: II**

(*на английском языке*)